\magnification 1200 \input miniltx \input graphicx.sty \makeatletter \def \Gin @driver{pdftex.def} \input color.sty \resetatcatcode \input amssym \input miniltx.tex
\font \bbfive = bbm5 \font \bbeight = bbm8 \font \bbten = bbm10     \font \eightbf = cmbx8 \font \eighti = cmmi8 \skewchar \eighti = '177
\font \eightit = cmti8 \font \eightrm = cmr8 \font \eightsl = cmsl8 \font \eightsy = cmsy8 \skewchar \eightsy = '60 \font \eighttt = cmtt8 \hyphenchar \eighttt = -1
  \font \sixi = cmmi6 \skewchar \sixi = '177 \font \sixrm = cmr6 \font \sixsy = cmsy6 \skewchar \sixsy = '60 \font \tensc =
cmcsc10   \scriptfont \bffam = \bbeight \scriptscriptfont \bffam = \bbfive \textfont \bffam = \bbten \newskip \ttglue
\def \eightpoint {\def \rm {\fam 0 \eightrm }\relax \textfont 0= \eightrm \scriptfont 0 = \sixrm \scriptscriptfont 0 = \fiverm \textfont 1 = \eighti \scriptfont 1 =
\sixi \scriptscriptfont 1 = \fivei \textfont 2 = \eightsy \scriptfont 2 = \sixsy \scriptscriptfont 2 = \fivesy \textfont 3 = \tenex \scriptfont 3 = \tenex
\scriptscriptfont 3 = \tenex \def \it {\fam \itfam \eightit }\relax \textfont \itfam = \eightit \def \sl {\fam \slfam \eightsl }\relax \textfont \slfam = \eightsl
\def \bf {\fam \bffam \eightbf }\relax \textfont \bffam = \bbeight \scriptfont \bffam = \bbfive \scriptscriptfont \bffam = \bbfive \def \tt {\fam \ttfam \eighttt
}\relax \textfont \ttfam = \eighttt \tt \ttglue = .5em plus.25em minus.15em \normalbaselineskip = 9pt \def \MF {{\manual opqr}\-{\manual stuq}}\relax \let \sc =
\sixrm \let \big = \eightbig \setbox \strutbox = \hbox {\vrule height7pt depth2pt width0pt}\relax \normalbaselines \rm } \def \TRUE {Y} \def \FALSE {N} \def \EMPTY
{} \def \ifundef #1{\expandafter \ifx \csname #1\endcsname \relax } \def \undefrule {\kern 2pt \vrule width 5pt height 5pt depth 0pt \kern 2pt} \def \UndefLabels {}
\def \possundef #1{\ifundef {#1}\undefrule {\eighttt #1}\undefrule \global \edef \UndefLabels {\UndefLabels #1\par } \else \csname #1\endcsname \fi } \newcount
\secno \secno = 0 \newcount \stno \stno = 0 \newcount \eqcntr \eqcntr = 0 \ifundef {showlabel} \global \def \showlabel {\FALSE } \fi \ifundef {auxwrite} \global
\def \auxwrite {\TRUE } \fi \ifundef {auxread} \global \def \auxread {\TRUE } \fi \def \define #1#2{\global \expandafter \edef \csname #1\endcsname {#2}} \long \def
\error #1{\medskip \noindent {\bf ******* #1}} \def \fatal #1{\error {#1\par Exiting...}\end } \def \advseqnumbering {\global \advance \stno by 1 \global \eqcntr
=0} \def \current {\ifnum \secno = 0 \number \stno \else \number \secno \ifnum \stno = 0 \else .\number \stno \fi \fi } \begingroup \catcode `\@=0 \catcode `\\=11
@global @def @textbackslash {\} @endgroup  \def \deflabel #1#2{\if \TRUE \showlabel \rem {\tiny [#1]}\fi \ifundef {#1PrimarilyDefined}\define
{#1}{#2}\define {#1PrimarilyDefined}{#2}\if \TRUE \auxwrite \immediate \write 1 {\textbackslash newlabel {#1}{#2}}\fi \else \edef \old {\csname #1\endcsname }\edef
\new {#2}\if \old \new \else \fatal {Duplicate definition for label ``{\tt #1}'', already defined as ``{\tt \old }''.}\fi \fi }  \def \label #1 {\deflabel {#1}{\current }} \def \equationmark #1 {\ifundef {InsideBlock} \advseqnumbering \eqno {(\current )} \deflabel {#1}{\current }
\else \global \advance \eqcntr by 1 \edef \subeqmarkaux {\current .\number \eqcntr } \eqno {(\subeqmarkaux )} \deflabel {#1}{\subeqmarkaux } \fi } \def \split
#1.#2.#3.#4;{\global \def \parone {#1}\global \def \partwo {#2}\global \def \parthree {#3}\global \def \parfour {#4}} \def \NA {NA} \def \ref #1{\split
#1.NA.NA.NA;(\possundef {\parone }\ifx \partwo \NA \else .\partwo \fi )}  \newcount \bibno \bibno = 0  \def \Bibitem #1 #2; #3; #4 \par {\smallbreak \global \advance \bibno by 1 \item {[\possundef {#1}]} #2, {``#3''}, #4.\par \ifundef
{#1PrimarilyDefined}\else \fatal {Duplicate definition for bibliography item ``{\tt #1}'', already defined in ``{\tt [\csname #1\endcsname ]}''.}  \fi \ifundef
{#1}\else \edef \prevNum {\csname #1\endcsname } \ifnum \bibno =\prevNum \else \error {Mismatch bibliography item ``{\tt #1}'', defined earlier (in aux file ?) as
``{\tt \prevNum }'' but should be ``{\tt \number \bibno }''.  Running again should fix this.}  \fi \fi \define {#1PrimarilyDefined}{#2}\if \TRUE \auxwrite
\immediate \write 1 {\textbackslash newbib {#1}{\number \bibno }}\fi } \def \jrn #1, #2 (#3), #4-#5;{{\sl #1}, {\bf #2} (#3), #4--#5} \def \Article #1 #2; #3; #4
\par {\Bibitem #1 #2; #3; \jrn #4; \par } \def \references {\begingroup \bigbreak \eightpoint \centerline {\tensc References} \nobreak \medskip \frenchspacing }
\catcode `\@=11 \def \c@itrk #1{{\bf \possundef {#1}}} \def \c@ite #1{{\rm [\c@itrk{#1}]}} \def \sc@ite [#1]#2{{\rm [\c@itrk{#2}\hskip 0.7pt:\hskip 2pt #1]}} \def
\du@lcite {\if \pe@k [\expandafter \sc@ite \else \expandafter \c@ite \fi } \def \cite {\futurelet \pe@k \du@lcite } \catcode `\@=12 \def \Headlines
#1#2{\nopagenumbers \headline {\ifnum \pageno = 1 \hfil \else \ifodd \pageno \tensc \hfil \lcase {#1} \hfil \folio \else \tensc \folio \hfil \lcase {#2} \hfil \fi
\fi }} \def \title #1{\medskip \centerline {\withfont {cmbx12}{\ucase {#1}}}}   \long \def \Quote #1\endQuote {\begingroup \leftskip 35pt \rightskip 35pt \parindent 17pt \eightpoint #1\par \endgroup } \long \def
\Abstract #1\endAbstract {\vskip 1cm \Quote \noindent #1\endQuote }   \def \Note #1{\footnote {}{\eightpoint #1}} \def \Date #1 {\Note {\it Date: #1.}}    \def \fix {\smallskip \noindent $\blacktriangleright $\kern 12pt}
\def \iskip {\medskip \noindent }    \def \ucase #1{\edef \auxvar {\uppercase
{#1}}\auxvar } \def \lcase #1{\edef \auxvar {\lowercase {#1}}\auxvar } \def \emph #1{{\it #1}} \def \section #1 \par {\global \advance \secno by 1 \stno = 0
\goodbreak \bigbreak \noindent {\bf \number \secno .\enspace #1.}  \nobreak \medskip \noindent } \def \state #1 #2\par {\begingroup \def \InsideBlock {} \medbreak
\noindent \advseqnumbering {\bf \current .\enspace #1.\enspace \sl #2\par }\medbreak \endgroup } \def \definition #1\par {\state Definition \rm #1\par } \newcount
\CloseProofFlag \def \closeProof {\eqno \endproofmarker \global \CloseProofFlag =1}  \long \def
\Proof #1\endProof {\begingroup \def \InsideBlock {} \global \CloseProofFlag =0 \medbreak \noindent {\it Proof.\enspace }#1 \ifnum \CloseProofFlag =0 \hfill
$\endproofmarker $ \looseness = -1 \fi \medbreak \endgroup } \def \quebra #1{#1 $$$$ #1} \def \explica #1#2{\mathrel {\buildrel \hbox {\sixrm #1} \over #2}} \def
\explain #1#2{\explica {\ref {#1}}{#2}} \def \Explain #1#2{\explain {#1}{#2} $$$$ \explain {#1}{#2}} \def \=#1{\explain {#1}{=}} \def \pilar #1{\vrule height #1
width 0pt} \def \stake #1{\vrule depth #1 width 0pt} \newcount \fnctr \fnctr = 0 \def \fn #1{\global \advance \fnctr by 1 \edef \footnumb {$^{\number \fnctr
}$}\footnote {\footnumb }{\eightpoint #1\par \vskip -10pt}} \def \text #1{\hbox {#1}}    \def \Item #1{\smallskip \item {{\rm #1}}} \newcount \zitemno \zitemno = 0 \def \izitem {\global \zitemno = 0}  \def \zitemplus {\global
\advance \zitemno by 1 \relax } \def \rzitem {\romannumeral \zitemno } \def \rzitemplus {\zitemplus \rzitem } \def \zitem {\Item {{\rm (\rzitemplus )}}}  \def \Zitem {\izitem \zitem } \def \zitemmark #1 {\deflabel {#1}{\current .\rzitem }{\def \showlabel {\FALSE }\deflabel {Local#1}{\rzitem }}} \def
\iItemmark #1 {\zitemmark {#1} } \def \zitemlbl #1 {\zitem \deflabel {#1}{\current .\rzitem }{\def \showlabel {\FALSE }\deflabel {Local#1}{\rzitem }}} \newcount
\nitemno \nitemno = 0  \def \nitem {\global \advance \nitemno by 1 \Item {{\rm (\number \nitemno )}}} \newcount \aitemno \aitemno = -1
\def \boxlet #1{\hbox to 6.5pt{\hfill #1\hfill }} \def \iaitem {\aitemno = -1} \def \aitemconv {\ifcase \aitemno a\or b\or c\or d\or e\or f\or g\or h\or i\or j\or
k\or l\or m\or n\or o\or p\or q\or r\or s\or t\or u\or v\or w\or x\or y\or z\else zzz\fi } \def \aitem {\global \advance \aitemno by 1\Item {(\boxlet \aitemconv )}}
\def \aitemmark #1 {\deflabel {#1}{\aitemconv }} \def \Bitem {\Item {$\bullet $}} \def \Case #1:{\medskip \noindent {\tensc Case #1:}} \def \<{\left \langle \vrule
width 0pt depth 0pt height 8pt } \def \>{\right \rangle } \def \({\big (} \def \){\big )} \def \ds {\displaystyle } \def \and {\hbox {\quad and \quad }} \def
\calcat #1{\,{\vrule height8pt depth4pt}_{\,#1}}  \def \imply {\mathrel {\Rightarrow }} \def \IFF {\kern
7pt\Leftrightarrow \kern 7pt} \def \IMPLY {\kern 7pt \Rightarrow \kern 7pt} \def \for #1{\quad \forall \,#1} \def \endproofmarker {\square } \def \"#1{{\it #1}\/}
 \def \inv {^{-1}} \def \*{\otimes } \def \caldef #1{\global \expandafter \edef \csname #1\endcsname {{\cal #1}}} \def \mathcal
#1{{\cal #1}} \def \bfdef #1{\global \expandafter \edef \csname #1\endcsname {{\bf #1}}} \bfdef N \bfdef Z \bfdef C \bfdef R  \def
\exists {\mathchar "0239\kern 1pt } \if \TRUE \auxread \IfFileExists {\jobname .aux}{\input \jobname .aux}{\null } \fi \if \TRUE \auxwrite \immediate \openout 1
\jobname .aux \fi \def \close {\if \EMPTY \UndefLabels \else \message {*** There were undefined labels ***} \iskip ****************** \ Undefined Labels: \tt \par
\UndefLabels \fi \if \TRUE \auxwrite \closeout 1 \fi \par \vfill \supereject \end } \def \bye {\close } \def \medsum {\mathop {\mathchoice {\hbox {$\mathchar
"1350$}}{\mathchar "1350}{\mathchar "1350}{\mathchar "1350}}} \def \medcup {\mathop {\mathchoice {\raise 1pt \hbox {$\mathchar "1353$}}{\mathchar "1353}{\mathchar
"1353}{\mathchar "1353}}} \def \medcap {\mathop {\mathchoice {\raise 1pt \hbox {$\mathchar "1354$}}{\mathchar "1354}{\mathchar "1354}{\mathchar "1354}}}     \def \newsection #1
\par {\global \edef \secname {#1}\global \advance \secno by 1 \stno = 0} \def \sectiontitle \par {\goodbreak \bigbreak \noindent {\bf \number \secno .\enspace
\secname .}  \nobreak \medskip \noindent } \input pictex \def \beginmypix #1, #2, #3, #4 /{ \begingroup \def \LeftLimit {#1} \def \RightLimit {#2} \def \TopLimit
{#3} \def \BotLimit {#4} \noindent \hfill \beginpicture \setcoordinatesystem units <0.0080truecm, -0.0080truecm> \setplotarea x from {\LeftLimit } to {\RightLimit
}, y from {\TopLimit } to {\BotLimit } \put {\null } at {\LeftLimit } {\TopLimit } \put {\null } at {\RightLimit } {\BotLimit } } \def \endmypix {\endpicture \hfill
\null \endgroup } \caldef G \def \restr #1{{\upharpoonright _{#1}}} \def \qi /{quasi-invariant} \def \dqi /{$D$-\qi /} \def \charac #1{1_{\kern -0.5pt #1}} \def \d
{\,d} \def \soma #1#2{\mathop {\hbox {${\mathchar "1358}$}}_{#1\in R(#2)}} \def \sminus {{\setminus }} \def \Gz {\G ^{(0)}} \def \Mplus {{\cal
M}^{^{\scriptscriptstyle +}}} \def \Mp {\Mplus \kern -2pt} \def \MpX {\Mp (X)} \def \MpXR {\Mplus _{R, \rho }(X)} \def \MpBX {\Mp \big (X,{\cal B}(X)\big )} \def
\prp /{proper} \def \gap /{{\tensc g\kern -0.5pt a\kern -0.5pt p}}


  \centerline {\bf QUASI-INVARIANT MEASURES FOR}
  \smallskip
  \centerline {\bf GENERALIZED APPROXIMATELY PROPER}
  \smallskip
  \centerline {\bf EQUIVALENCE RELATIONS}

  \bigskip
  \centerline {\tensc
  R. Bissacot\footnote {$^\ast $}{\eightrm Institute of Mathematics and
Statistics (IME-USP), University of S\~ao Paulo, Brazil.},
  R. Exel\footnote {$^{\ast \ast }$}{\eightrm Universidade Federal de Santa
Catarina and University of Nebraska. Partially supported by CNPq.},
  R. Frausino$^\ast $ and T. Raszeja$^\ast $}

\Abstract
  We introduce a generalization of the notion of approximately proper
equivalence relations studied by Renault and with it we build an \'etale
groupoid.  Choosing a suitable set of continuous functions to play the role
of a potential, we construct a cocycle in that groupoid and discuss the
corresponding Radon-Nikodym problem.
  \endAbstract

\section Introduction

In \cite {ApRenault}, Renault introduced the notion of an approximately
proper equivalence relation on a compact topological space $X$, consisting
of an increasing sequence $\{R_n\}_{n\in {\bf N}}$ of equivalence relations
on $X$, each of which is proper in the sense that the corresponding quotient
map is a local homeomorphism.  When equipped with the inductive limit
topology, the union $R=\bigcup _nR_n$ becomes an \'etale groupoid, and if
one is moreover given a suitable sequence of continuous real valued
functions on $X$, a cocycle\fn {In this paper the term \emph {cocycle} will
always be taken to mean a \emph {one}-cocycle.} may be defined on $R$.

As its title suggest, the main goal of \cite {ApRenault} is to study the
corresponding Radon-Nikodym problem, i.e., to find the probability measures
on $X$ which are \qi / with Radon-Nikodym derivative equal to the
aforementioned cocycle.

Among other things, the relevance os solving the Radon-Nikodym problem lies
in the fact that the solutions lead to KMS states on the groupoid C*-algebra
and hence have a profound relevance to Statistical Mechanics.

As mentioned in \cite [Section 7]{ApRenault}, approximately proper
equivalence relations arise naturally in the study of local homeomorphisms
from a compact topological space to itself.  Precisely speaking, given a
compact topological space $X$, and a local homeomorphism
  $
  \sigma :X\to X,
  $
  one lets,
  $$
  R_n = \big \{(x, y)\in X\times X: \sigma ^n(x)=\sigma ^n(y)\big \},
  $$
  for each $n\geq 0$,
  and it is not hard to see that each $R_n$ is a proper equivalence relation
so that $\{R_n\}_{n\in {\bf N}}$ is an approximately proper equivalence
relation in the sense of \cite {ApRenault}.

Proeminent examples of local homeomorphisms on compact topological spaces
are given by Markov shifts.  On the other hand, in a recent paper \cite
{PaperOne}, we have focused on a generalization of Markov shifts introduced
by M.~Laca and the second named author in \cite {infinoa}, which in turn
have been shown by Renault \cite {cuntzlike} to consist essentially of a
generalized \emph {shift space}, with the notable difference that the \emph
{shift map} is no longer defined on the whole space, but only on a proper
open subset.

The precise setup of \cite {cuntzlike} is that of a locally compact space
$X$, an open subset $U\subseteq X$, and a local homeomorphism
  $$
  \sigma :U\to X.
  $$
  However, if one starts from this data, it is not possible to build an
approximately proper equivalence relation by the procedure indicated above,
not least because $\sigma ^n$ fails to be defined on the whole space $X$.
If one wants to make sense of the relation
  $$
  x\sim y \iff \sigma ^n(x)=\sigma ^n(y),
  $$
  one must restrict attention to elements $x$ and $y$ for which $\sigma
^n(x)$ and $\sigma ^n(y)$ make sense, namely elements of the domain of
$\sigma ^n$, which we shal henceforth denote by $U_n$.  We thus define
  $$
  R_n = \big \{(x, y)\in U_n\times U_n: \sigma ^n(x)=\sigma ^n(y)\big \},
  $$
  which is clearly a proper equivalence relation on $U_n$.
  If $n\leq m$, observe that
  $$
  \sigma ^n(x)=\sigma ^n(y)\ \imply \ \sigma ^m(x)=\sigma ^m(y),
  $$
  as long as all of the above terms are defined, i.e, as long as $x$ and $y$
lie in the smaller set $U_m$.  This may be more sucintly expressed by saying
that
  $$
  R_n\cap (U_m\times U_m)\subseteq R_m.
  \equationmark Cresce
  $$

  If one misreads the above inclusion, ignoring the intersection with
$U_m\times U_m$, one will be left with the impression that the $R_n$ are
\emph {increasing}, just as in \cite {ApRenault}, although this is evidently
not true given that the sets where these relations are defined in fact \emph
{decrease}.

Another distinctive feature of the $R_n$ is the fact that, still under the
hypothesis that $n\leq m$, one has that $U_m$ is invariant under $R_n$,
meaning that
  $$
  \big ((x,y)\in R_n\big ) \ \wedge \ (y\in U_m)\ \imply \ x\in U_m,
  $$
  which may be expressed by saying that
  $$
  R_n\cap (U_n\times U_m) \subseteq U_m\times U_m.
  \equationmark Inva
  $$
  The reader may easily verify that, together, \ref {Cresce} and \ref {Inva}
are equivalent to
  $$
  R_n\cap (U_n\times U_m) \subseteq R_m,
  $$
  which might not have an immediately intuitive interpretation, but due to
its sheer simplicity, is adopted in this work as the main axiom in our
generalization of Renault's notion of approximately proper equivalence
relations, given in full detail in \ref {DefGap}, below, and referred to as
a \gap /, for short.

The main aim of the present work is to conduct a study of \gap /s along the
lines of Renault's study of approximately proper equivalence relations.  We
thus show that the union $R=\bigcup _nR_n$ is an equivalence relation, hence
a principal groupoid, which becomes \'etale when given the inductive limit
topology.  A suitable notion of potential is introduced, leading up to a
cocycle relative to which the Radon-Nikodym problem may be investigated.

Since each $R_n$ is assumed to be proper, one has that $R$ is the (not
necessarily increasing) union of proper equivalence relations, so it is not
surprising that the study of proper relations is as important here as it is
in \cite {ApRenault}.  Should our $U_n$ be compact we would be able to
borrow the results of the first few sections of \cite {ApRenault}, but the
example of infinite state Markov shifts, our main motivation, requires an
understanding of the Radon-Nikodym problem for proper relations on
non-compact spaces.  The lack of compactness indeed brings several
complications, most of them steming from the fact that equivalence classes
no longer need to be finite.  For example, the normalization achieved in
\cite [Proposition 3.1.iii]{ApRenault} by means of replacing a potential
$\rho '$ by
  $$
  \rho (x) := {\rho '(x) \over \sum _{y\sim x} \rho '(y)}
  $$
  needs to be dealt with in a more careful way if equivalence classes are
allowed to be infinite.

Once the proper case is taken care of, we apply our results for \gap /s,
showing, among other things, that \qi / measures may be characterized, much
in the same way as DLR measures, as those which are fixed by a family of
conditional expectations.  See section \ref {QIGapSect}, and in particular
Corollary \ref {MainGapQi}, for full details.

The existence part of the Radon-Nikodym problem, which follows easily from
compactness when that property is present, e.g.~as in \cite [8.2]{eq}, turns
out to be a delicate question here.  In fact existence may already fail in
the proper case, but it is nevertheless easy to determine precisely when
this happens.  The crucial point is to analyze the partition function
  $$
  \zeta (x) = \sum _{y\sim x} \rho (y),
  $$
  defined in \ref {DefinePsi}, which may well return infinite values, should
equivalence classes be infinite.  The set of points $x$ for which $\zeta
(x)=\infty $, which we denote by $Z_\rho $, is a forbidden zone for finite
\qi / measures in the sense that any such measure assigns zero mass to
$Z_\rho $.  Thus, if $\zeta $ is identically infinite, a situation very easy
to arrange, there are no nontrivial solutions for the Radon-Nikodym problem.
Excluding this extreme situation, i.e.~when $\zeta $ is finite on at least
one point, one may easily show the existence of \qi / measures.  See section
\ref {QionEtale} for more details in the proper case.

Unfortunately we have no definitive answer for the existence question in the
most general situation of \gap /s treated here, which is perhaps to be
expected given that similar results rely heavily on compactness.  However we
can offer several partial existence results which the reader will find in
section \ref {ExistenceSection}.

The second named author would like to acknowledge the warm hospitality of
Rodrigo Bissacot and his group during a very productive visit to the
University of S\~ao Paulo, when the bulk of the results presently being
reported were developped.

\newsection Proper equivalence relations

\sectiontitle

As mentioned above, we start by analyzing proper equivalence relations,
avoiding the compactness assumption.

\state {Standing Hypothesis} \label Standing
  \rm Throughout this notes we will assume that $X$ is a locally compact,
second countable, metrizable space.

We will denote the $\sigma $-algebra of Borel measurable subset of $X$ by
${\cal B}(X)$, and the set of all Borel measurable functions
  $$
  f:X\to [0,+\infty ]
  $$ by $\MpBX $.

\definition
  An equivalence relation $R\subseteq X\times X$ is said to be \emph {\prp
/}, provided the quotient space $X/R$ is Hausdorff, and the quotient map
  $$
  \pi :X\to X/R
  $$
  is a local homeomorphism\fn {A map $\varphi :X\to Y$, between topological
spaces $X$ and $Y$, is said to be a \emph {local homeomorphism} provided for
every $x$ in $X$, there are open subsets $A\subseteq X$, and $B\subseteq Y$,
such that $x\in A$ and $\varphi $ is a homeomorphism from $A$ onto $B$.}.

\fix From now on we shall fix a \prp / equivalence relation $R$ on $X$.

\medskip Given any $x$ in $X$, we will denote its equivalence class by
$R(x)$, in symbols
  $$
  R(x) = \{y\in X:(x, y)\in R\}.
  $$

\definition \label IntroE
  For each $f$ in $\MpBX $, we will let
  $$
  E(f)\calcat x = \medsum _{y\in R(x)} f(y).
  $$

Observe that the above sum could very well diverge, in which case we of
course set $E(f)\calcat x\stake {6pt}$ to be $\infty $. Therefore, like $f$,
one has that $E(f)$ is a function taking values in $[0,\infty ]$.

We will soon prove that $E(f)$ is ${\cal B}(X)$-measurable, but so far we
will see it simply as an element of $\Mp \big (X,{\cal P}(X)\big )$, where
${\cal P}(X)$ is the $\sigma $-algebra of \emph {all} subsets of $X$ (with
respect to which any function is measurable).  In other words, $E$ may be
seen as a map
  $$
  E:\MpBX \to \Mp \big (X,{\cal P}(X)\big ).
  $$

\state Proposition \label SigmaAdditOfE
  $E$ is \emph {$\sigma $-additive}, in the sense that it is positively
homogeneous and
  $$
  E \Big (\medsum _{n=1}^\infty f_n \Big ) = \medsum _{n=1}^\infty E (f_n),
  $$
  for any sequence $\{f_n\}_n$ in $\MpBX $.

\Proof It is evident that $E$ is positively homogeneous.  Given any sequence
$\{f_n\}_n$ in $\MpBX $, for every $x$ in $X$, we have
  $$
  E \big (\medsum _{n=1}^\infty f_n\big )\calcat x=
  \medsum _{y\in R(x)} \medsum _{n=1}^\infty f_n(y) =
  \medsum _{n=1}^\infty \ \medsum _{y\in R(x)} f_n(y) =
  \medsum _{n=1}^\infty E(f_n)\calcat x.
  \closeProof
  $$
  \endProof

\state Proposition \label EPreservrMeasblty
  If $f$ is in $\MpBX $, then so is $E(f)$.

\Proof Let $\{U_n\}_{n\in {\bf N}}$ be a countable open cover of $X$, such
that the quotient map
  $$
  \pi :X\to X/R
  $$
  is a homeomorphism when restricted to each $U_n$.  Also let $\{\psi
_n\}_{n\in {\bf N}}$ be a partition of unit subordinate to this cover.

Given $f$ in $\MpBX $, put $f_n=f\psi _n$, so that $f=\sum _nf_n$,
pointwise.  Using \ref {SigmaAdditOfE}, we then have that
  $$
  E(f) = E \Big (\medsum _{n=1}^\infty f_n\Big ) = \medsum _{n=1}^\infty
E(f_n),
  $$
  so it suffices to prove that each $E(f_n)$ is Borel-measurable.  Write
  $$
  \tau :\pi (U_n)\to U_n
  $$
  for the inverse of the restriction of $\pi $ to $U_n$, and let $V_n=\pi
\inv \big (\pi (U_n)\big )$.  We then claim that
  $$
  E(f_n)\calcat x = \left \{\matrix {
  f_n\big (\tau (\pi (x))\big ), & \hbox { if } x\in V_n, \hfill \cr \pilar
{12pt}
  0, & \hbox { otherwise.}
  }\right .
  $$

  Indeed, when $x$ is not in $V_n$, then $\pi (x)$ is not in $\pi (U_n)$.
So, while $f_n$ vanishes outside $U_n$, there is no $y$ in $U_n$ such that
$(x,y)\in R$.  The sum defining $E(f_n)\calcat x$ therefore has no nonzero
terms, and hence $E (f_n)\calcat x=0$.

On the other hand, if $x$ is in $V_n$, then $\pi (x)\in \pi (U_n)$, so $\pi
(x)=\pi (y)$, for a unique $y$ in $U_n$, namely $y=\tau \big (\pi (x)\big
)$, whence $f(y)$ is the only possibly nonzero term in the aforementioned
sum.  Therefore
  $$
  E(f_n)\calcat x = f(y) = f_n\big (\tau (\pi (x))\big ),
  $$
  as claimed.
  Since the correspondence $x\mapsto f_n\big (\tau (\pi (x))\big )$ is
easily seen to be Borel-measura\-ble on $V_n$, we have that $E(f_n)$ is
Borel-measurable on $X$.
  \endProof

Notice that, in view of the above result, $E $ may be viewed as a map from
$\MpBX $ to itself.  We therefore no longer need to consider the $\sigma
$-algebra ${\cal P}(X)$, and we shall henceforth use the simplified notation
  $$
  \MpX := \MpBX .
  $$

A few other useful properties of $E $ are as follows:

\state Proposition \label UsefulProp
  Given $f,g\in \MpX $, one has that:
  \izitem
  \zitem $E(f)$ is $R$-invariant, meaning that if $(x,y)\in R$, then
$E(f)\calcat x = E(f)\calcat y$,
  \zitem if $g$ is $R$-invariant, then $E(gf)=gE(f)$,
  \zitemlbl ERCresc if $f\leq g$, then $E(f)\leq E(g)$,
  \zitemlbl EROrb if $f$ vanishes outside a subset $A\subseteq X$, then
$E(f)$ vanishes outside
  $$
  \hbox {Orb}(A) := \{y\in X: \exists x\in A,\ (x, y)\in R\}.
  $$

\Proof We prove only \ref {LocalEROrb}.
  Given $x\in X$, suppose that
  $$
  0\neq E(f)\calcat x = \medsum _{y\in R(x)} f(y).
  $$
  Then it is easy to see that there exists at least one $y$ such that
$(x,y)\in R$, and $f(y)\neq 0$.  Consequently $y\in A$, whence $x\in \hbox
{Orb}(A)$.  This proves that
  $$
  E(f)\calcat x\neq 0\imply x\in \hbox {Orb}(A),
  $$
  from where the conclusion follows immediately.
  \endProof

When multiplying extended real numbers, as in the multiplication ``$gf$''
above, we adopt the convention according to which $0\times \infty =\infty
\times 0=0$.  A trivial, but highly relevant fact to be noted regarding this
convention is that multiplication of positive extended real numbers is both
associative and infinitely distributive, i.e.,
  $$
  c\medsum _{n=1}^\infty a_n = \medsum _{n=1}^\infty ca_n,
  $$
  for every $c$ and every sequence $\{a_n\}_n$ in $[0,\infty
]$. Incidentally the above choice for the value of $0\times \infty $ is
necessary for the validity of the distributive property, since
  $$
  0\times \infty = 0\times \medsum _{n=1}^\infty {1\over n} = \medsum
_{n=1}^\infty 0\times {1\over n} = 0.
  $$

Given $\rho , f$ in $\MpX $, we have that $E(\rho f)\in \MpX $.  Fixing
$\rho $ we may then define the map
  $$
  E_\rho :f\in \MpX \mapsto E(\rho f)\in \MpX ,
  $$
  which is clearly also $\sigma $-additive.  Therefore, for every measure
$\nu $ on ${\cal B}(X)$, we may consider the measure $E_\rho ^*(\nu )$ given
by \ref {Daniell}.  Some elementary observations regarding $E_\rho ^*(\nu )$
are in order:

\state Proposition \label NewDuasIntegrais
  Given a function $\rho $ in $\MpX $ as well as a measure $\nu $ on ${\cal
B}(X)$, one has that
  \izitem
  \zitem $E_\rho ^*(\nu )$ is a finite measure if and only if $E(\rho )$ is
$\nu $-integrable,
  \zitem if $A$ is any
  $R$-invariant\fn {A subset $A\subseteq X$ is said to be \emph
{$R$-invariant} if, whenever $x\in A$, and $(x, y)\in R$, one has that $y\in
A$.}
  Borel-measurable subset of $X$ with $\nu (A)=0$, then $E_\rho ^*(\nu
)(A)=0$ as well.

\Proof The first point follows immediately from
  $$
  E_\rho ^*(\nu )(X)=\int _X 1\d E_\rho ^*(\nu ) = \int _X E_\rho (1)\d \nu
= \int _X E(\rho )\d \nu .
  $$

  Regarding (ii), and denoting the characteristic function of $A$ by
$\charac A$, it is obvious that $\rho \charac A$ vanishes outside
$A$. Therefore $E(\rho \charac A)$ vanishes outside $\hbox {Orb}(A)$ by \ref
{EROrb}. However, since $A$ is invariant, we have that $\hbox {Orb}(A) = A$,
so in fact $E(\rho \charac A)$ vanishes outside $A$.  Therefore
  $$
  E_\rho ^*(\nu )(A) =
  \int _X \charac A \d E_\rho ^*(\nu ) =
  \int _X E(\rho \charac A) \d \nu =
  \int _{X\sminus A} E(\rho \charac A) \d \nu =
  0.
  \closeProof
  $$
  \endProof

\newsection The operator $E $ on $C_c(X)$

\sectiontitle

In this section we continue assuming that $X$ satisfies \ref {Standing} and
that $R\subseteq X\times X$ is a \prp / equivalence relation on $X$.
Whenever we speak of $R$ as a topological space, we will be referring to the
topology induced on $R$ by the product topology of $X\times X$.

We will often view $R$ as a \emph {groupoid} under the multiplication
operation according to which the product $(x,y)\cdot (z,w)$ is defined if
and only if $y=z$, in which case it is set to be $(x,w)$.  The unit space of
such a groupoid is therefore the diagonal $\{(x, x): x\in X\}$, which we
identify with $X$ in the obvious way.  The \emph {range} and \emph {source}
maps are then given respectively by
  $$
  r(x,y)=x, \and s(x,y)=y, \for (x,y)\in R.
  $$
  It is well known that $R$ is then a Hausdorff \'etale groupoid.

\state Proposition \label Continuity
  Given any continuous, complex valued function $f$ on $R$, suppose that $f$
vanishes outside a given subset $L\subseteq R$, such that $s(L)$ is
relatively compact.  Then:
  \izitem
  \zitem
  The expression
  $$
  g(y) = \sum _{\gamma :r(\gamma )=y} f(\gamma )
  $$
  gives a well defined and continuous function on $X$.
  \zitem If $r(L)$ is also relatively compact, then $g$ has compact support.

\Proof We first claim that $R$ is closed in $X\times X$.  In order to see
this let
  $$
  \pi :X\to X/R
  $$
  denote the quotient map and observe that
  $$
  R =
  \{(x, y)\in X\times X: \pi (x)=\pi (y)\} =
  \{(x, y)\in X\times X: \big (\pi (x), \pi (y)\big )\in \Delta \},
  $$
  where $\Delta $ is the diagonal in $X/R\times X/R$.  Since $X/R$ is
Hausdorff, we have that $\Delta $ is closed, and hence $R$ is closed in
$X\times X$.

  Given $y_0$ in $X$, let $K$ be a compact neighborhood of $y_0$, and
observe that
  $$
  r\inv (K)\cap L\subseteq K\times s(L),
  $$
  so $r\inv (K)\cap L$ is relatively compact in $X\times X$, and hence also
in the closed subspace $R$.  For each
  $$
  \gamma \in \overline {r\inv (K)\cap L},
  $$
  namely the closure of $r\inv (K)\cap L$ within $R$, let $U_\gamma
\subseteq R$ and $V_\gamma \subseteq X$ be open sets such that $\gamma \in
U_\gamma $, $r(\gamma )\in V_\gamma $, and such that the restriction of the
range map $r$ to $U_\gamma $ gives a homeomorphism onto $V_\gamma $.

We shall also insist that, whenever $r(\gamma )\neq y_0$, the open
neighborhood $V_\gamma $ of $r(\gamma )$ is chosen such that $y_0\notin
\overline {V}_\gamma $.  We therefore get an open cover $\{U_\gamma
\}_\gamma $ of the compact set $\overline {r\inv (K)\cap L}$, from which one
may extract a finite subcover, say $\{U_\gamma \}_{\gamma \in F}$, where $F$
is some finite set of $\gamma $'s.

Splitting $F$ according to whether or not $r(\gamma )=y_0$, we define
  $$
  F_1 = \{\alpha \in F: r(\alpha )=y_0\}, \and F_2 = \{\beta \in F: r(\beta
)\neq y_0\}.
  $$ We then put
  $$
  V = \hbox {int}(K)\cap \bigcap _ {\alpha \in F_1} V_\alpha \cap \bigcap _
{\beta \in F_2} X\setminus \overline {V}_\beta ,
  $$
  and we claim that $y_0\in V$.  To see this, notice that $y_0$ lies in
$\hbox {int}(K)$ because $K$ is a neighborhood of $y_0$.  Moreover, for
every $\alpha $ in $F_1$, we have that
  $$
  y_0=r(\alpha ) \in V_\alpha ,
  $$
  and finally, for every $\beta $ in $F_2$, we have explicitly chosen
$V_\beta $ so that $y_0\notin \overline {V}_\beta $.

We next claim that, for every $\eta $ in $R$,
  $$
  r(\eta )\in V \ \wedge \ f(\eta )\neq 0 \IMPLY \eta \in U_\alpha , \hbox {
for some } \alpha \in F_1.
  \equationmark Implicacao
  $$ Indeed, given $\eta $ satisfying the above antecedent, we clearly have
that
  $$
  \eta \in r\inv (K)\cap L,
  $$
  so there exists some $\gamma $ in $F$, such that $\eta \in U_\gamma $, and
we would now like to decide whether $\gamma $ lies in $F_1$ or in $F_2$.
The key observation here is that
  $$
  r(\eta )\in V\cap r(U_\gamma ) = V\cap V_\gamma ,
  $$
  so $V_\gamma $ has a nonempty intersection with $V$, and this can only
happen when $\gamma \in F_1$, thus completing the proof of \ref
{Implicacao}.

Next use the fact that $R$ is Hausdorff to produce a collection of pairwise
disjoint open sets $\{W_\alpha \}_{\alpha \in F_1}$ such that each $\alpha
\in W_\alpha $ and finally put
  $$
  \Omega =V\cap \bigcap _{\alpha \in F_1} r(W_\alpha \cap U_\alpha ).
  $$

  Noticing that $W_\alpha \cap U_\alpha $ is open in $U_\alpha $, we see
that $r(W_\alpha \cap U_\alpha )$ is open in $r(U_\alpha )=V_\alpha $, so
$\Omega $ is an open subset of $X$.  Also, since $\alpha \in W_\alpha \cap
U_\alpha $, we have that $y_0=r(\alpha )\in r(W_\alpha \cap U_\alpha )$, so
$y_0\in \Omega $.

  For each $\alpha $ in $F_1$, denote by $t_\alpha $ the inverse of the
homeomorphism
  $$
  r\restr {U_\sigma }:U_\alpha \to V_\alpha ,
  $$
  and, regarding the function $g$ referred to in the statement,
 we claim that for every $y$ in $\Omega $, one has that
  $$
  g(y) = \sum _{\alpha \in F_1} f\big (t_\alpha (y)\big ).
  \equationmark NewG
  $$
  To prove this claim it is enough to show that
  $$
  \{\gamma \in R: r(\gamma )=y,\ f(\gamma )\neq 0\} = \{t_\alpha (y): \alpha
\in F_1,\ f(t_\alpha (y))\neq 0\},
  \equationmark SetsInSum
  $$
  and that the $t_\alpha (y)$ in the description of the set in the right
hand side above are pairwise distinct.

With respect to this last statement, notice that for each $\alpha $ in
$F_1$,
  $$
  t_\alpha (y) \in t_\alpha \big (r(W_\alpha \cap U_\alpha )\big ) =
W_\alpha \cap U_\alpha \subseteq W_\alpha ,
  $$
  so the $t_\alpha (y)$ lie in pairwise disjoint sets and hence are
necessarily pairwise distinct.

We next observe that the inclusion ``$\supseteq $'' in \ref {SetsInSum} is
evident, so we focus on the reverse inclusion ``$\subseteq $''.  For this,
pick $\gamma $ in $R$ such that $r(\gamma )=y$, and $f(\gamma )\neq 0$, and
notice that by \ref {Implicacao} it follows that $\gamma \in U_\alpha $ for
some $\alpha \in F_1$.  Therefore
  $\gamma =t_\alpha (y)$, so we see that $\gamma $ lies in the set in the
right hand side of \ref {SetsInSum}.

This proves \ref {SetsInSum} and hence also \ref {NewG}, from where it is
clear that the sum defining $g$ has finitely many nonzero terms, so that $g$
is well defined, and moreover that
 $g$ is continuous.

In order to prove (ii), it is enough to observe that if $g(y)\neq 0$, there
must be at least one $\gamma $ with $r(\gamma )=y$, and $f(\gamma )\neq 0$,
whence $\gamma \in L$, and then $y\in r(L)$. Viewing through the
counterpositive
  $$
  y\notin r(L) \IMPLY f(\gamma )=0, \hbox { for all $\gamma $ such that
$r(\gamma )=y$},
  $$
  which in turn implies that $g(y)=0$.  Thus $g$ vanishes outside the
relatively compact set $r(L)$, and hence it is compactly supported.
  \endProof

\state Corollary \label SumSource
  Given $f\in C_c(R)$, the correspondence
  $$
  x\mapsto \sum _{y\in R(x)} f(x,y)
  $$
  defines a compactly supported, continuous function on $X$.

\Proof Follows immediately from \ref {Continuity.ii} upon choosing $L$ to be
the support of $f$.
  \endProof

Recall that the operator $E$ defined in \ref {IntroE} is only defined for
non-negative functions.  This is due to the fact that the summation involved
in its definition is not supposed to converge but, as long as the summands
are non-negative, we may always assign a sensible value to the sum, that
value being $\infty $ in the divergent case.  In case of compactly supported
functions the situation is however much better behaved:

\state Proposition \label IntroECont
  Given $f$ in $C_c(X)$, and for every $x$ in $X$, the sum
  $$
  \sum _{y\in R(x)} f(y),
  $$
  has at most finitely many nonzero terms.  Moreover, defining
  $$
  E(f)\calcat x = \sum _{y\in R(x)} f(y), \for x\in X,
  $$
  one has that $E(f)$ is a continuous function on $X$.

\Proof Since $R$ is a \prp / equivalence relation, we have that $R(x)$ is a
closed, discrete set for every $x$ in $X$.  Therefore, if $K$ is the compact
support of $f$, one has that $R(x)\cap K$ is finite from where the first
assertion follows immediately.

Addressing the last assertion, consider the continuous function
  $$
  g:(x, y)\in R\mapsto f(y)\in {\bf C}.
  $$
  Denoting the support of $f$ by $L$, notice that $g$ vanishes outside the
set
  $$
  K:=(X\times L)\cap R.
  $$

  Since $s(K)\subseteq L$, we see that $s(K)$ is relatively compact. We may
therefore employ \ref {Continuity.i} to conclude that the function $g$
defined there is continuous, namely
  $$
  g(x) =
  \sum _{\gamma :r(\gamma )=x} g(\gamma ) =
  \sum _{y\in R(x)} g(x, y) =
  \sum _{y\in R(x)} f(y) =
  E(f)\calcat x,
  $$
  concluding the proof.
  \endProof

In view of the above result we get a map
  $$
  E :C_c(X) \to C(X).
  $$
  On the other hand, recall that in \ref {IntroE} we defined an operator
  $$
  E :\MpX \to \MpX ,
  $$
  using the exact same formula as in \ref {IntroECont}.  Clearly the two
operators referred to above coincide on the intersection of their domains,
so there is no ambiguity in using the same notation ``$E $'' for these maps.

Some of the main properties of $E $ on $C_c(X)$ reflect those listed in \ref
{UsefulProp}:

\state Proposition \label PropERho
  For every $f$ and g in $C_c(X)$, one has that:
  \izitem
  \zitemlbl ERInvar $E(f)$ is $R$-invariant,
  \zitemlbl ERCondExp if $g$ is $R$-invariant, then $E(gf)=gE(f)$,
  \zitemlbl ERCont $E(f)$ is continuous,
  \zitemlbl ERLimit $E(f)$ is bounded.

\Proof Leaving the easy proofs of (i) and (ii) to the reader, we notice that
\ref {LocalERCont} was already proved in \ref {IntroECont}.

Regarding \ref {LocalERLimit}, let $K$ be the compact support of $f$, and
let $M$ be the supremum of $|E(f)|$ on $K$, which is finite by \ref
{LocalERCont}.  We will then prove that $|E(f)|$ is bounded by $M$ on all of
$X$.  In order to prove that
  $$
  \big |E(f)\calcat x\big | \leq M,
  \equationmark Objetivo
  $$
  for any given $x$ in $X$, we may evidently assume that $E(f)\calcat x\neq
0$.  In this case
  $$
  0 \neq E(f)\calcat x = \sum _{y\in R(x)} f(y),
  $$
  so there exists at least one $y$ in $K$ such that $(x, y)\in R$.
Therefore
  $$
  \big |E(f)\calcat x\big | \explain {LocalERInvar}= \big |E(f)\calcat
{y}\big | \leq M,
  $$
  proving \ref {Objetivo}.
  \endProof

Complementing \ref {NewDuasIntegrais}, we may now describe a few other
relevant properties of $E_\rho ^*(\nu )$ under the extra hypothesis that
$\rho $ is finitely valued and continuous.

\state Proposition \label MoreOnMuNu
  Let $\nu $ be a measure on ${\cal B}(X)$, and let
  $
  \rho :X\to {\bf R},
  $
  be a non-negative, continuous function. Setting $\mu =E_\rho ^*(\nu )$,
one has that:
  \izitem
  \zitem if $\nu (X)<\infty $, then $\mu $ is a Borel measure (i.e.~finite
on compact sets),
  \zitem if $\mu '$ is any measure on ${\cal B}(X)$ such that
  $$
  \infty >\int _Xf\d \mu ' = \int _XE(\rho f)\d \nu , \for f\in C_c^+(X),
  $$
  then $\mu =\mu '$, and in particular the above identity holds for every
$f$ in $\MpX $.

\Proof In order to verify (i), and using \ref {ConditionBorel}, it is enough
to prove that every $f$ in $C_c^+(X)$ is $\mu $-integrable.  Given such an
$f$, notice that the continuity of $\rho $ implies that $\rho f$ lies in
$C_c(X)$, whence $E(\rho f)$ is bounded by \ref {ERLimit}.  Therefore
  $$
  \int _X f\d \mu =\int _XE(\rho f)\d \nu \leq \nu (X)\Vert E(\rho f)\Vert
_\infty < \infty .
  $$

Addressing (ii), observe that the hypothesis says that every $f$ in
$C_c^+(X)$ is $\mu '$-integra\-ble, so another application of \ref
{ConditionBorel} tells us that $\mu '$ is a Borel measure.  By hypothesis we
then have that
  $$
  \int _X f\d \mu ' =\int _X f\d \mu , \for f\in C_c^+(X),
  $$
  from where we deduce that $\mu $ is also a Borel measure.  Since any $f$
in $C_c(X)$ may be written as the linear combination of functions in
$C_c^+(X)$, we deduce that the identity displayed above holds for every $f$
in $C_c(X)$, so $\mu =\mu '$, by \ref {AllRegular} and the uniqueness part
of the Riesz-Markov Theorem.
  \endProof

\newsection Proper equivalence relations and \qi / measures

\def \zi {\zeta ^{\scriptscriptstyle -1}}

\sectiontitle

  \label QionEtale
  \fix As before, throughout this section we fix a space $X$ satisfying \ref
{Standing}, as well as a \prp / equivalence relation $R$ on $X$.  We will
moreover fix a continuous function
  $$
  \rho :X\to {\bf R},
  $$
  which will henceforth be supposed \emph {strictly positive}, i.e,
  $$
  \rho (x)>0, \for x\in X,
  $$
  and which will be referred to as the \emph {potential} for $R$.

The relevance of $\rho $ is that it leads to a multiplicative cocycle on $R$
via the formula
  $$
  D(x, y):=\rho (x)/\rho (y), \for (x,y)\in R,
  \equationmark IntroCocycle
  $$
  and the goal of this section is to study \emph {\qi /} measures relative
to this cocycle.  See \ref {DefQuaseInvariant} below for the precise
definition.

\medskip In some applications of our theory, the role of $\rho $ is played
by the function $\rho (x) = e^{\beta h(x)}$, where $\beta >0$ and $h$ is a
continuous, real valued function on $X$.  The assumption that $\rho $ is
strictly positive then holds automatically.  Another reason why we need to
assume that $\rho $ is never zero is that, otherwise, the the above
definition of $D$ would run into trouble.

\definition \label DefinePsi
  For each $x$ em $X$, define
  $$
  \zeta (x) := E(\rho )\calcat x = \soma yx \rho (y).
  $$
  We will refer to $\zeta $ as the \emph {partition function} for the
potencial $\rho $.

\state Proposition \label PsiLimitInfer
  $\zeta $ is bounded below by $\rho $, and consequently
  $$
  0<\zeta (x)\leq \infty , \for x\in X.
  $$

\Proof Obvious.
  \endProof

Since $\zeta $ is defined to be $E(\rho )$, we have by \ref
{EPreservrMeasblty} that $\zeta $ lies in $\MpX $.  We may in fact prove
that $\zeta $ satisfies a stronger regularity property:

\state Proposition \label LSC
  $\zeta $ is lower semi-continuous.

\Proof Let $\{\varphi _n\}_n$ be as in \ref {ApproxUnit}.  Then
  $$
  \zeta =
  E(\rho ) =
  E \big (\lim _{n\to \infty } \rho \varphi _n\big ) \={MonotoneForE}
  \lim _{n\to \infty } E(\rho \varphi _n),
  $$
  whence $\zeta $ is the limit of an increasing sequence of continuous
functions by \ref {IntroECont}, from where the conclusion follows.
  \endProof

\state Corollary \label IntroZRho
  The set
  $$
  Z_\rho = \{x\in X: \zeta (x) = \infty \}
  $$
  is a $ G_\delta $, hence a Borel set.

\Proof Noting that
  $$
  Z_\rho = \bigcap _{k\in {\bf N}} \{x\in X :\zeta (x)>k\},
  $$
  the result is an immediate consequence of \ref {LSC}.
  \endProof

In what follows we will make frequent references to the function $\zi $, so
it is worth discussing it briefly now.  Recall from \ref {PsiLimitInfer}
that $\zeta (x)>0$, for all $x$ in $X$, so we will never run into the
trouble of considering the inverse of zero.  On the other hand, when $\zeta
(x)=\infty $, we evidently put $\zi (x)=0$.

In view of our convention that $\infty \times 0=0$, observe that
  $$
  \zeta (x)\zeta (x)\inv = \left \{\matrix {
  0, & \hbox {if } x\in Z_\rho , \hfill \cr
  1, & \hbox {otherwise.}
  } \right .
  \equationmark ZeroInfinito
  $$
  so we have that
  $$
  \zeta \zi = \charac {X\sminus Z_\rho }.
  \equationmark PropYr
  $$

In particular we note the following partial-isometric-like property of
$\zeta $, to be used shortly:
  $$
  \zi \zeta \zi =\zi .
  \equationmark PartialIso
  $$

  All things considered, we will see that the somewhat unusual fact that
$\zeta \zi $ vanishes on $Z_\rho $ will not be so crucial.  For example, we
will soon encounter expressions such as
  $$
  \int _X \zeta \zi \d \mu ,
  $$
  but often the measure $\mu $ will also vanishes on $Z_\rho $, so the funny
behavior of $\zeta \zi $ on $Z_\rho $ becomes totally irrelevant.

We next recall the definition of a \qi / measure in the special case of
\'etale groupoids.

\definition \label DefQuaseInvariant
  \cite [I.3.15]{Renault}
  Let $\G $ be an \'etale groupoid and let $D:\G \to \R _+^*$ be a
multiplicative cocycle.  A measure $\mu $ on $\Gz $ is said to be \emph {\qi
/} relative to $D$ when
  $$
  \int _{\Gz } \sum _ {\gamma \in r\inv (x)} f(\gamma ) \d \mu (x)= \int
_{\Gz } \sum _ {\gamma \in s\inv (x)} f(\gamma ) D(\gamma ) \d \mu (x),
  $$
  for every $f$ in $C_c(\G )$.

For the case of our groupoid $R$, the above quasi-invariance condition
becomes
  $$
  \int _X \soma yx f(x,y) \d \mu (x)= \int _X \soma xy f(x,y) D(x,y) \d \mu
(y),
  \equationmark OriginalQuaseInvariant
  $$
  for every $f$ in $C_c(R)$.

\medskip The following result lists several equivalent conditions for a
measure to solve the Radon-Nikodym problem.

\state Theorem \label MainResult
  Let $X$ be a topological space satisfying \ref {Standing}. Also let $R$ be
a \prp / equivalence relation on $X$, seen as an \'etale groupoid.  Given a
continuous, strictly positive function $\rho : X \to {\bf R}$, consider the
cocycle defined on $R$ by $D(x,y)=\rho (x)/\rho (y)$.  Then, for every
finite measure $\mu $ on $X$, the following are equivalent:
  \def \tbox #1#2#3{\hbox to 2.1cm{\hfill $\ds #1$\hfill }$\ =\ $\hbox to
2.7cm{$\ds #2, $}$\forall #3$
    \vrule depth 15pt width 0pt}
  \izitem
  \zitem $\mu $ is \dqi /,
  \zitem \tbox {\int _X fE (\rho g)\d \mu }{\int _X E (\rho f)g\d \mu }{f,
g\in C_c(X)},
  \zitemlbl Harmonic \tbox {\int _X f\zeta \d \mu }{\int _X E (\rho f) \d
\mu }{f\in C_c^+(X)},
  \zitemlbl Idempot \tbox {\int _X f\d \mu }{\int _X E (f\rho \zi )\d \mu
}{f\in C_c^+(X)},
  \zitemlbl Built there exists a positive measure $\nu $ on $X$, with
respect to which $\zeta $ is integrable, and \hfill \break \tbox {\int _X
f\d \mu }{\int _X E (\rho f)\d \nu }{f\in C_c(X)}.
  \medskip \noindent In addition, if any of the above equivalent conditions
hold, then $\mu (Z_\rho )=0$.

\Proof (i) $\Rightarrow $ (ii).  Pick $f$ and $g$ in $C_c(X)$, and consider
the function $F$ on $R$ given by the formula
  $$
  F(x,y)=f(x)g(y)\rho (y).
  $$
  Plugging $F$ in \ref {OriginalQuaseInvariant} we have
  $$
  \int _X \soma yx f(x)g(y)\rho (y) \d \mu (x)= \int _X \soma xy f(x)g(y)
\rho (x) \d \mu (y),
  $$
  which translated precisely into (ii).

\medskip \noindent (ii) $\Rightarrow $ \ref {LocalHarmonic}.  Given $f$ in
$C_c^+(X)$, let $\{\varphi _n\}_n$ be as in \ref {ApproxUnit}.  Then
  $$
  \int _X f\zeta \d \mu =
  \int _X fE(\rho ) \d \mu =
  \int _X fE(\lim _{n\to \infty } \rho \varphi _n) \d \mu \= {MonotoneForE}
$$$$=
  \lim _{n\to \infty } \int _X fE(\rho \varphi _n) \d \mu \explica {(ii)}=
  \lim _{n\to \infty } \int _X E(\rho f)\varphi _n \d \mu =
  \int _X E(\rho f) \d \mu .
  $$

\medskip \noindent \ref {LocalHarmonic} $\Rightarrow $ \ref {LocalIdempot}.
  We will first prove that $\mu (Z_\rho )=0$.  In order to do it suppose by
way of contradiction that $\mu (Z_\rho )>0$.  Since $\mu $ is finite, it is
regular by \ref {AllRegular}, so there exists a compact set $K\subseteq
Z_\rho $, such that $\mu (K)>0$.  Using Uryhson, take $f$ in $C_c^+(X)$,
such that $f|_K=1$, so that
  $$
  \infty = \infty \times \mu (K) =
  \int _K \zeta \d \mu =
  \int _K f\zeta \d \mu \quebra \leq
  \int _X f\zeta \d \mu \explica {\ref {LocalHarmonic}}=
  \int _X E(\rho f) \d \mu \leq
  \Vert E(\rho f)\Vert _\infty \ \mu (X) \explain {ERLimit}< \infty .
  $$
  Arriving at a contradiction we conclude that $\mu (Z_\rho )=0$, as
desired.

We next claim that \ref {LocalHarmonic} indeed holds for all $f$ in $\MpX $,
namely that
  $$
  \int _X f\zeta \d \mu = \int _X E(\rho f) \d \mu , \for f\in \MpX .
  \equationmark ExtendedHarmonc
  $$

Letting $\zeta \mu $ and $\mu $ play the roles of $\mu '$ and $\nu $,
respectively, in \ref {MoreOnMuNu.ii}, we only need to prove that every $f$
in $C_c^+(X)$ is integrable with respect to $\zeta \mu $, but this follows
from
  $$
  \int _Xf\zeta \d \mu \explica {\ref {LocalHarmonic}}=
  \int _XE_{\rho } (f)\d \mu \leq \mu (X)\Vert E_\rho (f)\Vert _\infty
\explain {ERLimit}<\infty .
  $$
  Therefore \ref {ExtendedHarmonc} is verified so, for any
  $f$ in $C_c^+(X)$, we may plug in $f\zi $ there, obtaining
  $$
  \int _XE(\rho f\zi )\d \mu =
  \int _Xf\zi \zeta \d \mu =
  \int _Xf\d \mu ,
  $$
  where the last equality is justified by \ref {ZeroInfinito}, which says
that $\zeta \zi =1$ on $X\setminus Z_\rho $, and by the fact that $\mu $
vanishes on $Z_\rho $.
  This proves \ref {LocalIdempot}.

\medskip \noindent \ref {LocalIdempot} $\Rightarrow $ \ref {LocalBuilt}.
  Defining $\nu :=\zi \mu $, and given $f$ in $C_c^+(X)$, we have that
  $$
  \int _X f\d \mu =
  \int _X E(\rho f\zi )\d \mu \={UsefulProp}
  \int _X E(\rho f)\zi \d \mu = \int _X E(\rho f)\d \nu .
  $$
  Since $C_c(X)$ is linearly spanned by $C_c^+(X)$, the last assertion in
(v) follows.  Furthermore, employing \ref {MoreOnMuNu.ii} once more, one has
that
  $$
  \int _X f\d \mu = \int _X E(\rho f)\d \nu , \for f\in \MpX ,
  $$
  so we are allowed to plug $f=1$ above, whence
  $$
  \int _X\zeta \d \nu = \int _XE(\rho )\d \nu = \int _X 1 \d \mu = \mu
(X)<\infty ,
  $$
  hence proving the remaining first assertion of \ref {LocalBuilt}.

\medskip \noindent \ref {LocalBuilt} $\Rightarrow $ (i).  In order to prove
that $\mu $ is \dqi /, we need to check \ref {OriginalQuaseInvariant} for
every $f$ in $C_c(R)$.  As a notational aid, let us temporarily write
  $$
  A(x) = \soma yx f(x,y), \and B(y) = \soma xy f(x,y) D(x,y),
  $$
  so that our goal is to prove that $A$ and $B$ have the same integral
relative to $\mu $.  En passant, notice that $A$ and $B$ lie in $C_c(X)$ by
\ref {SumSource}.

Starting from the left-hand-side of \ref {OriginalQuaseInvariant}, observe
that
  $$
  \int _X \soma yx f(x,y) \d \mu (x)=
  \int _X A \d \mu \={LocalBuilt}
  \int _X E(\rho A) \d \nu \quebra =
  \int _X \soma zx\rho (z)\soma yz f(z,y) \d \nu (x)=
  \int _X \soma zx\ \soma yz \rho (z)f(z,y) \d \nu (x).
  \equationmark LadoEsqFinal
  $$
  On the other hand, starting from the right-hand-side of \ref
{OriginalQuaseInvariant}, we have
  $$
  \int _X \soma zy f(z,y) D(z,y) \d \mu (y) =
  \int _X B \d \mu \explica {\ref {LocalBuilt}}=
  \int _X E(\rho B) \d \nu \quebra =
  \int _X \soma yx\rho (y)\soma zy f(z,y) D(z,y)\d \nu (x) =
  \int _X \soma yx\ \soma zy \rho (z)f(z,y) \d \nu (x).
  \equationmark LadoDirFinal
  $$

Notice that the diference between \ref {LadoEsqFinal} and \ref
{LadoDirFinal} is simply that, in the former, the sum ranges over all pairs
$(z,y)$ such that
  $$
  x\mathop {\sim }_Rz\mathop {\sim }_Ry,
  $$
  while, in the latter, the pairs $(y,z)$ considered are those for which
  $$
  x\mathop {\sim }_Ry\mathop {\sim }_Rz.
  $$

  Being an equivalence relation, $R$ is transitive, whence in both cases
above the sum ranges over all $y$ and all $z$ in the equivalence class of
$x$, and therefore we see that \ref {LadoEsqFinal} and \ref {LadoDirFinal}
coincide.  This proves \ref {OriginalQuaseInvariant} and hence that $\mu $
is \dqi /.
  \endProof

The characterization given by \ref {Built} may be used to produce \dqi /
measures, as we now show:

\state Corollary \label ExistenciaQInvar
  Given a measure $\nu $ on $X$ such that $\zeta $ is $\nu $-integrable,
there exists a unique finite, \dqi / measure $\mu $ on $X$ such that
  $$
  \int _Xf\d \mu = \int _XE(\rho f)\d \nu , \for f\in C_c(X).
  $$

\Proof Given $\nu $, let $\mu =E_\rho ^*(\nu )$, so
  $$
  \mu (X) = \int _X 1 \d \mu = \int _XE(\rho )\d \nu = \int _X\zeta \d \nu
<\infty ,
  $$
  so $\mu $ is indeed a finite measure and it is \dqi / because it satisfies
\ref {Built}.  The uniqueness of $\mu $ now follows from the uniqueness part
of the Riesz-Markov Theorem.
  \endProof

The next result settles the question regarding the existence of nontrivial
\dqi / measures.

\state Corollary \label CondParaExist
  The following are equivalent:
  \izitem
  \zitem there exists at least one \dqi / probability measure on $X$,
  \zitem $\zeta $ is not identically infinite.

\Proof Recall that $Z_\rho $ is the set of points where $\zeta $ is
infinite, so (ii) is equivalent to saying that $Z_\rho \neq X$, or
equivalently that $X\sminus Z_\rho $ is nonempty.

  Assuming (i), let $\mu $ be a \dqi / probability measure on $X$.  By the
last sentence in \ref {MainResult} we have that $\mu (Z_\rho )=0$, and hence
that $\mu (X\sminus Z_\rho )=1$, so $X\sminus Z_\rho \neq \emptyset $,
proving (ii).  Conversely, if $X\sminus Z_\rho $ is nonempty, it is easy to
exhibit a measure $\nu $ on $X$ satisfying
  $$
  \int _X \zeta \d \nu =1.
  $$
  Take, for example, any point $y_0\in X\sminus Z_\rho $ and, observing that
  $
  0<\zeta (y_0)<\infty
  $
  by \ref {PsiLimitInfer}, it is enough to choose
  $$
  \nu =\zeta (y_0)\inv \delta _{y_0},
  $$
  where $\delta _{y_0}$ is the Dirac measure on $y_0$.
  Given any such $\nu $, the measure $\mu $ built in \ref {ExistenciaQInvar}
in terms of $\nu $ is a \dqi / probability measure, proving (i).
  \endProof

The remainder of this section will be devoted to a closer look at the fourth
condition of \ref {MainResult}.

\state Proposition \label Projector
  Consider the operator
  $
  P_\rho :\MpX \to \MpX ,
  $
  given by
  $$
  P_\rho (f) = E(f\rho \zi ), \for f\in \MpX .
  $$
  Then
  \izitem
  \zitem $P_\rho (1) = \charac {X\sminus Z_\rho }$,
  \zitem $ P_\rho ^2=P_\rho $,
  \zitem the range of $P_\rho $ coincides with the set $\MpXR $, consisting
of all $R$-invariant functions $f$ in $\MpX $ which vanish on $Z_\rho $.

\Proof We should first observe that, since $\rho $ and $\zi $ lie in $\MpX
$, the range of $P_\rho $ is indeed a subset of $\MpX $ by \ref
{EPreservrMeasblty}.

In order to prove the first assertion, we compute
  $$
  P_\rho (1) = E(\rho \zi ) \={UsefulProp.ii}
  E(\rho )\zi =
  \zeta \zi \={ZeroInfinito}
  \charac {X\sminus Z_\rho }.
  $$ We next claim that:
  \iaitem
  \aitem the range of $P_\rho $ is \emph {contained} in $\MpXR $, and
  \aitem $P_\rho (f)=f$, for every $f$ in $\MpXR $.

  \smallskip
  In order to verify (a), pick any $f$ in $\MpX $.  Since
  $$
  P_\rho (f) = E(f\rho )\zi
  $$
  by \ref {UsefulProp.ii}, and since $\zi $ vanishes on $Z_\rho $, then
$P_\rho (f)$ also vanishes on $Z_\rho $.  The fact that $P_\rho (f)$ lies in
$\MpXR $ then follows immediately from \ref {UsefulProp.i}.

To prove (b), let $f\in \MpXR $.  Then
  $$
  P_\rho (f) =
  E(f\rho \zi ) \={UsefulProp.ii}
  fE(\rho )\zi =
  f\zeta \zi =
  f,
  $$
  where the last step relies on the fact that $f$ vanishes on $Z_\rho $.
  This said, (ii) and (iii) follow trivially from (a) and (b).
  \endProof

Among the characterizations of \dqi / measures given by \ref {MainResult}, a
particularly useful one is \ref {MainResult.iv}, given the nice properties
of the operator $P_\rho $ described in \ref {Projector}.  For that reason,
and also for future reference, we restate part of the conclusions of \ref
{MainResult} in a way as to emphasize the importance of $P_\rho $.

\state Corollary \label MainResultWithP
  Under the conditions of \ref {MainResult} one has that $\mu $ is \dqi / if
and only if $P_\rho ^*(\mu )=\mu $.

Some further important facts involving $P_\rho ^*$ are as follows.

\state Proposition \label PropOfPstar
  Let $\nu $ be any finite measure on $X$.  Then
  \izitem
  \zitem if $A\subseteq X$ is an invariant Borel subset, then $P_\rho ^*(\nu
)(A) =\nu (A\sminus Z_\rho )$,
  \zitem if $\nu $ vanishes on an invariant Borel set $A\subseteq X$, then
the same is true for $P_\rho ^*(\nu )$,
  \zitem $P_\rho ^*(\nu )$ is finite,
  \zitem $P_\rho ^*(\nu )$ is nonzero if and only if $\nu (X\sminus Z_\rho
)$ is nonzero,
  \zitem if $\nu $ is a probability measure vanishing on $Z_\rho $, then so
is $P_\rho ^*(\nu )$,
  \zitem $P_\rho ^*\big (P_\rho ^*(\nu )\big )=P_\rho ^*(\nu )$,
  \zitem $P_\rho ^*(\nu )$ is \dqi /.

\Proof Given $A$ as in (i), we have
  $$
  P_\rho ^*(\nu )(A) =
  \int _X\charac A\d P_\rho ^*(\nu ) =
  \int _XE(\charac A\rho \zi )\d \nu \={UsefulProp.ii}
  \int _X\charac AE(\rho )\zi \d \nu \quebra =
  \int _X\charac A\zeta \zi \d \nu \={ZeroInfinito}
  \int _X\charac A\charac {X\sminus Z_\rho }\d \nu = \nu (A\sminus Z_\rho ).
  $$
  Points (ii--v) then follow immediately from (i).  Regarding (vi), it is an
easy consequence of \ref {Projector.ii}.  Finally let us prove (vii).  For
this, set $\mu =P_\rho ^*(\nu )$, so we see from (vi) that $P_\rho ^*(\mu
)=\mu $, and the conclusion follows from \ref {MainResultWithP}.
  \endProof

\newsection Generalized approximately proper equivalence relations

\sectiontitle

As before, throughout this section we assume that $X$ is a locally compact,
second countable, metrizable space.

\definition \label DefGap
  By a \emph {generalized approximately proper equivalence relation} on $X$,
a \gap / for short, we shall mean a pair
  $$
  {\cal R} = \Big (\{U_n\}_{n\in {\bf N}}, \{R_n\}_{n\in {\bf N}}\Big ),
  $$
  where each $U_n$ is an open subset of $X$, and each $R_n$ is a \prp /
equivalence relation on $U_n$, such that
  \izitem
  \zitem $X=U_0\supseteq U_1\supseteq U_2\supseteq \cdots $
  \zitem $R_0$ is the identity relation on $U_0$, that is, $R_0$ is the
diagonal in $U_0\times U_0$,
  \zitemlbl Crazy if $n\leq m$, then $R_n\cap (U_n\times U_m) \subseteq
R_m$.

Two immediate consequences of the definition are as follows:

\state Proposition \label AlternativeCrazy
  If ${\cal R}= \big (\{U_n\}_{n\in {\bf N}}, \ \{R_n\}_{n\in {\bf N}}\big
)$ is a \gap / on $X$ then, whenever $n\leq m$, one has that
  \izitem
  \zitemlbl Increasing the restriction of $R_n$ to $U_m$, namely $R_n\cap
(U_m\times U_m)$, is contained in $R_m$,
  \zitemlbl UBigInvar If $n\leq m$, then $U_m$ is \emph {invariant} under
$R_n$ in the sense that if a point $x$ in $U_n$ is equivalent under $R_n$ to
a point $y$ in $U_m$, then $x$ lies in $U_m$.

\Proof We have
  $$
  R_n\cap (U_m\times U_m) \subseteq R_n\cap (U_n\times U_m) \explain
{Crazy}\subseteq R_m,
  $$
  proving (i).
  If $x$ and $y$ are as in (ii), then
  $$
  (x, y)\in R_n\cap (U_n\times U_m) \explain {Crazy}\subseteq R_m \subseteq
U_m\times U_m,
  $$
  so $x\in U_m$.
  \endProof

It is not hard to see that also \ref {AlternativeCrazy.i--ii} imply \ref
{Crazy}, so the reader might think of the latter as subsuming the former,
which some may consider a more natural set of conditions.

The main motivation and the main source of examples for \gap /s is described
in detail in Section \ref {ExamplesSection}, below.

\fix From now on we fix a \gap / ${\cal R}= \big (\{U_n\}_{n\in {\bf N}}, \
\{R_n\}_{n\in {\bf N}}\big )$ on $X$.

\state Proposition
  Setting
  $$
  R = \medcup _{n\in {\bf N}}R_n,
  $$
  one has that $R$ is an equivalence relation on $X$.

\Proof The only slightly nontrivial point regards the transitivity of $R$.
In order to prove it, suppose that $(x,y)$ and $(y,z)$ lie in $R$.  We may
then pick $n$ and $m$ such that $(x,y)\in R_n$ and $(y,z)\in R_m$, and we
may assume without loss of generality that $n\leq m$.  In that case we have
that
  $$
  (x,y)\in R_n\cap (U_n\times U_m) \explain {Crazy}\subseteq R_m.
  $$
  Since $R_m$ is transitive we have that
  $$
  (x,z)\in R_m\subseteq R.
  \closeProof
  $$
  \endProof

\state Lemma \label OpenRel
  Equipping each $R_n$ with the topology induced from the product topology
on $X\times X$, one has that $R_n\cap R_m$ is open in $R_n$, for all $n$ and
$m$ in ${\bf N}$.

\Proof Given $(x,y)$ in $R_n\cap R_m$, by the definition of the product
topology on $R_n$ we must prove the existence of open subsets $V, W\subseteq
X$, such that
  $$
  (x,y)\in R_n\cap (V\times W) \subseteq R_n\cap R_m.
  $$

Assuming that $n\leq m$, choose $V=U_n$, and $W=U_m$, and observe that the
above inclusion is then immediately verified thanks to \ref {Crazy}.  On the
other hand, proving the result under the opposite assumption, i.e.~that
$n\geq m$, is equivalent to maintaining the assumption that $n\leq m$ (with
which the reader must be used to by now) and proving instead that
  $$
  (x,y)\in R_m\cap (V\times W) \subseteq R_n\cap R_m.
  \equationmark ProveInstead
  $$

  For each $k\in {\bf N}$, denote by $\pi _k$ the quotient map
  $$
  \pi _k:U_k\to U_k/R_k,
  $$
  and for each $k\in \{n, m\}$, let us choose an open set $W_k\subseteq
U_k$, such that $y\in W_k$, and such that $\pi _k$ restricts to a
homeomorphism from $W_k$ to the open set $\pi _k(W_k)$.  Replacing both
$W_n$ and $W_m$ by
  $$
  W:= W_n\cap W_m,
  $$
  we may assume that $W_n=W_m$.

  Notice that $\pi _n(x) = \pi _n(y) \in \pi _n(W)$, so we have that $x\in
\pi _n\inv \big (\pi _n(W)\big )$, and upon setting
  $$
  V:= \pi _n\inv \big (\pi _n(W)\big )\cap U_m,
  $$
  we see that $V$ is an open subset of $X$, and we moreover claim that \ref
{ProveInstead} holds.  The first part, namely that $(x,y)\in R_m\cap
(V\times W)$, is evident and, in order to prove that
  $$
  R_m\cap (V\times W) \subseteq R_n\cap R_m,
  \equationmark InsteadInclusion
  $$
  let us pick $(z,w)$ in the set appearing in the left-hand side above.  It
follows that $z\in V$, whence $\pi _n(z)\in \pi _n(W)$, so there exists some
$w'$ in $W$ such that $\pi _n(z)= \pi _n(w')$.  Another way to express this
is by saying that $(z,w')\in R_n$, but since $(z,w')$ also lies in
$U_m\times U_m$, we deduce that
  $$
  (z,w')\in R_n\cap (U_m\times U_m)\subseteq R_m.
  $$
  Recall that $(z,w)\in R_m$, as well, so transitivity yields $(w,w')\in
R_m$.  Observing that both $w$ and $w'$ lie in $W$, and using that $\pi _m$
is injective on $W$, we see that $w=w'$, whence
  $$
  (z,w) = (z,w') \in R_n.
  $$
  This finishes the verification of \ref {InsteadInclusion}, and hence also
of \ref {ProveInstead}, concluding the proof.
  \endProof

Recal that the inductive limit topology on the union of an increasing
sequence of topological spaces
  $$
  X_0\subseteq X_1\subseteq X_2\subseteq \cdots
  $$
  is the topology according to which a subset $U\subseteq \bigcup _nX_n$ is
open if and only if $U\cap X_n$ is open in $X_n$, for every $n$.  In our
situation, where $R=\bigcup _nR_n$, the $R_n$ do not form an increasing
sequence of subsets but one may nevertheless equipp $X$ with the topology
defined as above, that is, in which a subset $U\subseteq R$ is open if and
only if $U\cap R_n$ is open in $R_n$, for every $n$.  Even though this might
constitute a slight abuse of the language, we shall refer to that topology
as the \emph {inductive limit topology} on $R$.

\state Lemma
  Equipping $R$ with the inductive limit topology we have that each $R_n$ is
open in $R$.

\Proof Follows immediately from \ref {OpenRel}.
  \endProof

The two previous Lemmas form the key to showing the following result, whose
otherwise easy proof we leave for the reader.

\state Proposition
  Given a generalized approximately proper equivalence relation on $X$, one
has that $R$ is an \'etale groupoid when equipped with the inductive limit
topology.

\newsection Quasi-invariant measures and \gap /s

\def \Ninf {{\bf N}\cup \{\infty \}} \def \infunion {\kern -12pt\bigsqcup
_{\kern 12pt n\in \Ninf } \kern -12pt}

\sectiontitle

\label QIGapSect As before, throughout this section we assume that $X$ is a
locally compact, second countable, metrizable space.  We will also assume
that we are given a \gap /
  $$
  {\cal R} = \Big (\{U_n\}_{n\in {\bf N}}, \{R_n\}_{n\in {\bf N}}\Big )
  $$
  on $X$.

\definition \label DefinePotential
  By a \emph {potential}\/ for $\cal R$ we shall mean a collection
$\{k_n\}_{n\geq 1}$, of continuous functions
  $$
  k_n:U_n\to {\bf R},
  $$
  such that for every $n\geq 1$, one has that
  $$
  (x,y)\in R_{n-1}\cap (U_n\times U_n)\imply k_n(x)=k_n(y).
  \equationmark PotInv
  $$

\bigskip It is perhaps worth pointing out that a potential involves no
$k_0$.  On the other hand, the lowest case of \ref {PotInv} is tautological,
that is, when $n=1$, we have that $R_{n-1}$, also known as $R_0$, is the
identity relation, and it is no surprise that $x=y$ implies that
$k_n(x)=k_n(y)$.

Regarding \ref {PotInv}, notice that
  $$
  R_{n-1}\cap (U_n\times U_n) = R_{n-1}\cap R_n,
  $$
  because
  $$
  R_{n-1}\cap R_n \subseteq R_{n-1}\cap (U_n\times U_n) \subseteq
R_{n-1}\cap (U_{n-1}\times U_n) \explain {Crazy}\subseteq R_{n-1}\cap R_n.
  $$ An equivalent way to state \ref {PotInv} is therefore to require that
$k_n(x)=k_n(y)$, for all
  $(x,y)$ in $R_{n-1}\cap R_n.$

\fix From now on we assume that we are given a potential $\{k_n\}_{n\geq 1}$
for $\cal R$.

\medskip Our next goal is to use a potential to produce a cocycle on the
groupoid
  $
  R=\medcup _{n\in {\bf N}} R_n.
  $
  As a first step we introduce the following notation:

\definition
  For all $n\geq 0$, let $h_n:U_n\to {\bf R}$, be defined recursively by
$h_0=0$, and
  $$
  h_n=h_{n-1}\restr {U_n}+k_n, \for n\geq 1.
  $$
  In addition, for all $n\geq 0$, we define
  $$
  c_n:(x,y)\in R_n\mapsto h_n(x)-h_n(y)\in {\bf R}.
  $$

\bigskip Of course one may alternatively define $h_n$ by
  $$
  h_n=\sum _{i=1}^n k_i\restr {U_n},
  $$
  observing that, when $n=0$, the usual convention about sums without any
summands gives $h_0=0$, as expected.

\state Proposition \label OneStepCn
  For every $n\geq 1$, and all $(x,y)$ in $R_{n-1}\cap R_n$, one has that
  $$
  c_{n-1}(x,y) = c_n(x,y).
  $$

\Proof The difference between $c_n(x,y)$ and $c_{n-1}(x,y)$ is precisely
$k_n(x)-k_n(y)$, but since
  $(x,y)$ lies in $R_{n-1}\cap R_n$, condition \ref {PotInv} applies.
  \endProof

\state Proposition \label GapCocycle
  There exists a (necessarily unique) continuous cocycle $c$ on $R$, such
that $c=c_n$ on each $R_n$.

\Proof We first claim that, whenever $0\leq n\leq m$, one has that
  $$
  R_n\cap R_m = R_n\cap R_{n+1}\cap \ldots \cap R_{m-1}\cap R_m.
  \equationmark RChain
  $$
  In order to see this, it is clearly enough to show that $R_n\cap
R_m\subseteq R_k$, for every $k$ with $n\leq k\leq m$, and in turn this
follows from
  $$
  R_n\cap R_m \subseteq
  R_n\cap (U_m\times U_m) \subseteq R_n\cap (U_n\times U_k) \explain
{Crazy}\subseteq R_k.
  $$

Given any $\gamma $ in $R$, choose $n$ such that $\gamma \in R_n$, and put
  $$
  c(\gamma )=c_n(\gamma ).
  $$
  To see that this is well defined, suppose that $\gamma \in R_m$, for some
other $m$, and let us prove that $c_n(\gamma )=c_m(\gamma )$.  Assuming
without loss of generality that $n\leq m$, it follows from \ref {RChain}
that $\gamma \in R_n\cap R_{n+1}\cap \ldots \cap R_{m-1}\cap R_m$, so we may
apply \ref {OneStepCn} to show that
  $$
  c_n(\gamma )=c_{n+1}(\gamma )=\cdots =c_{m-1}(\gamma )=c_m(\gamma ).
  $$

  This shows that $c$ is well defined and we leave it as an easy exercise to
show that $c$ is a continuous cocycle on $R$.
  \endProof

We shall next present two general results about \qi / measures on \'etale
groupoids, to be used later.

\state Proposition \label UnionSubGpds
  Let $\G $ be an \'etale groupoid and suppose that we are given a
collection $\{\G _i\}_{i\in I}$ of open subgroupoids $\G _i\subseteq \G $,
such that $\G =\bigcup _{i\in I}\G _i$.  Suppose moreover that $D:\G \to
{\bf R}^*_+$ is a continuous multiplicative cocycle and that $\mu $ is a
finite measure on $\G ^{(0)}$.  Then $\mu $ is \dqi / if and only if the
restriction (see \ref {RestrictingMeasures} for a discussion regarding the
concept of restricting a measure to a subset) of $\mu $ to $\G _i^{(0)}$ is
\qi / relative to the restriction of $D$ to $\G _i$, for every $i$.

\Proof We prove only the ``if'' part, leaving the ``only if'' part to the
reader.  We must therefore check \ref {DefQuaseInvariant} for every $f$ in
$C_c(\G )$.  By \cite [3.10]{actions} (which holds even if $\G $ is
non-Hausdorff), we have that $f$ may be written as a finite linear
combination of functions $f_j$, each of which lies in $C_c(U_j)$, for some
open bissection $U_j$.
  Therefore, since both sides of \ref {DefQuaseInvariant} are clearly linear
with respect to $f$, it suffices to prove \ref {DefQuaseInvariant} under the
assumption that $f\in C_c(U)$, for some open bissection $U$.

Letting $K$ be the compact support of $f$, recall that the hypotheses imply
that $\{U\cap \G _i\}_i$ is an open cover for $K$.  Choosing a finite
subcover $\{U\cap \G _{i_k}\}_{k=1}^n$ and a partition of unit $\{\varphi
_k\}_{k=1}^n$ subordinate to it \cite [21.1.5]{RK}, we may write
  $$
  f=\sum _{k=1}^nf\varphi _k,
  $$
  observing that $f\varphi _k$ lies in $C_c(\G _{i_k})$.

The upshot of this argument is that we may further reduce \ref
{DefQuaseInvariant} by assuming that $f$ is supported on a single $\G _i$.
  Under this assumption, observe that the integrands in both sides of \ref
{DefQuaseInvariant} vanish whenever $x$ is not in $\G _i^{(0)}$, so it
suffices to verify a variant of \ref {DefQuaseInvariant}, namely where both
occurences of $\G ^{(0)}$ are replaced by $\G _i^{(0)}$.  The resulting
expression is then seen to hold because the restriction of $\mu $ to $\G
_i^{(0)}$ is \dqi / by hypothesis.
  \endProof

\state Lemma \label InvarSetAndMeasure
  Let $\G $ be an \'etale groupoid with a continuous multiplicative cocycle
$D:\G \to {\bf R}^*_+$, and let $\mu $ be a finite, \dqi / measure on $\Gz
$.
  \izitem
  \zitem If $\varphi $ is a bounded,
  invariant\fn {A function $\varphi $ defined on $\Gz $ is said to be
invariant when $\varphi (r(\gamma ))=\varphi (s(\gamma ))$, for every
$\gamma $ in $\G $.},
  Borel-measurable function on $\Gz $, then $\varphi \mu $ is also \dqi /.
  \zitem If $E$ is an invariant\fn {A subset $E\subseteq \Gz $ is said to be
invariant when $r(\gamma )\in E\Leftrightarrow s(\gamma )\in E$, for every
$\gamma $ in $\G $.}, Borel subset of $\Gz $, then $\mu _E := \charac E\,
\mu $ is also \dqi /.

\Proof In order to prove (i) we pick any $f$ in $C_c(\G )$ and we set out to
verify \ref {DefQuaseInvariant} relative to the measure $\varphi \mu $.
Starting from the left-hand side, we have
  $$
  \int _{\Gz } \sum _ {\gamma \in r\inv (x)} f(\gamma ) \d \varphi \mu (x)=
  \int _{\Gz } \sum _ {\gamma \in r\inv (x)} f(\gamma )\varphi (x) \d \mu
(x)\quebra =
  \int _{\Gz } \sum _ {\gamma \in r\inv (x)} f(\gamma )\varphi \big
(r(\gamma )\big ) \d \mu (x) \={DefQuaseInvariant}
  \int _{\Gz } \sum _ {\gamma \in s\inv (x)} f(\gamma )\varphi \big
(s(\gamma )\big ) D(\gamma ) \d \mu (x) \quebra =
  \int _{\Gz } \sum _ {\gamma \in s\inv (x)} f(\gamma ) D(\gamma ) \varphi
(x)\d \mu (x) =
  \int _{\Gz } \sum _ {\gamma \in s\inv (x)} f(\gamma ) D(\gamma ) \d
\varphi \mu (x),
  $$
  proving (i).
  Point (ii) now follows from (i) upon taking $\varphi $ to be the
characteristic function of $E$.
  \endProof

Returning to the \gap / we have fixed at the beginning of this section, and
assuming we are given a potential $\{k_n\}_{n\geq 1}$, leading up to the
cocycle $c$ of \ref {GapCocycle}, consider the multiplicative cocycle
  $$
  D:\gamma \in R\mapsto e^{c(\gamma )}\in {\bf R}^*_+,
  $$
  as well as the multiplicative cocycles
  $$
  D_n:\gamma \in R_n\mapsto e^{c_n(\gamma )}\in {\bf R}^*_+.
  $$
  We then have the following immediate consequence of \ref {UnionSubGpds}:

\state Corollary \label ReduceToEtale
  Let $\mu $ be a finite measure on $X$. Then $\mu $ is \dqi / if and only
if $\mu \restr {U_n}$ is $D_n$-\qi / for every $n$ in ${\bf N}$.

Once the question of the quasi-invariance of a measure $\mu $ on $X$ is
reduced to the quasi-invariance of measures on \prp / equivalence relations,
namely the $R_n$'s in the above Corollary, the results of Section \ref
{QionEtale} apply.
  Our goal in what follows is to patch the conclusions of these results for
the various $R_n$ in a meaningful way from the point of view of $R$.
  For each $n\in {\bf N}$, we shall let
  $$
  \rho _n:x\in U_n\mapsto e^{h_n(x)}\in {\bf R}_+^*,
  $$
  and we will henceforth let
  $\zeta _n$ be the partition function given by \ref {DefinePsi} in terms of
$\rho _n$.

Alongside $\rho _n$, $D_n$ and $\zeta _n$, all of the other ingredients
introduced in Section \ref {QionEtale} will also be relevant here, such as
$Z_{\rho _n}$ and $P_{\rho _n}$, as well as the operator $E_n$ on $\Mp
(U_n)$ given by \ref {IntroE} relative to the equivalence relation $R_n$.

As a first use of these notations we have the following immediate
consequence of \ref {ReduceToEtale} and \ref {MainResultWithP}.

\state Corollary \label MainGapQi
  Let $\mu $ be a finite measure on $X$. Then $\mu $ is \dqi / if and only
if
  $$
  P_{\rho _n}^*( \mu \restr {U_n} )= \mu \restr {U_n} , \for n\in {\bf N}.
  $$

Part of the difficulty in simultaneously dealing with so many maps and sets
is the fact that they each refer to a different equivalence relation.
Attempting to bring everything to a common environment we introduce the
following:

\definition \label ExtendedDefs
  Let $n\in {\bf N}$ be given.
  \izitem
  \zitem For any $f$ in $\Mp (U_n)$, we will denote by $\iota _n(f)$ the
extension of $f$ to the whole of $X$ obtained by setting it to be zero
outside $U_n$.  When no confusion is likely to arise, we shall denote that
extension simply by $f$, by abuse of language.
  \zitem We will write $Z_n$ and $Y_n$ for $Z_{\rho _n}$ and $U_n\sminus
Z_{\rho _n}$, respectively, and we will view both $Z_n$ and $Y_n$ as subsets
of $X$ (which of course they are).
  \zitem We will denote by $F_n$ the map from $\MpX $ to itself given, for
every $f$ in $\MpX $, and for every $x$ in $X$, by
  $$
  F_n(f)\calcat x = \left \{\matrix {
  \ds \medsum _{y\in R_n(x)} f(y), & \hbox {if } x\in U_n, \phantom {\rho
_n(y)\zeta _n(y)\inv }\cr \pilar {12pt}
  0, & \hbox {otherwise.}\hfill
  }\right .
  $$
  \zitem We will denote by $Q_n$ the map from $\MpX $ to itself given, for
every $f$ in $\MpX $, and for every $x$ in $X$, by
  $$
  Q_n(f)\calcat x = \left \{\matrix {
  \ds \medsum _{y\in R_n(x)} f(y)\rho _n(y)\zeta _n(y)\inv , & \hbox {if }
x\in U_n, \cr \pilar {12pt}
  0, & \hbox {otherwise.}
  }\right .
  $$
  \zitemlbl PartialInvar We will say that a given $f$ in $\MpX $ is \emph
{$R_n$-invariant}\/ if $f(x)=f(y)$, whenever $(x,y)\in R_n$.

\state Remarks \label RemarksOnExtendedDefs
  \rm
  \iaitem
  \aitem Since $h_0=0$, we have that $\rho _0=1$, and clearly also $\zeta
_0=1$.  Therefore $Z_0=\emptyset $.
  \aitem Notice that $F_n(f)$ and $Q_n(f)$ may be alternatively defined as
  $$
  F_n(f) = \iota _n\big (E_n(f\restr {U_n})\big ),\and Q_n(f) = \iota
_n(P_{\rho _n}\big (f\restr {U_n})\big ).
  $$
  For that reason $F_n$ and $Q_n$ should be seen as natural extensions of
$E_n$ and $P_{\rho _n}$ to $\MpX $, respectively.  Notice also that
  $$
  Q_n(f) = F_n\big (f\iota _n(\rho _n\zi _n)\big ), \for f\in \MpX .
  $$
  \aitem Observe that the invariance of a function under an equivalence
relation is a concept usually considered when the relation is defined on the
\emph {whole domain} of said function.  However, the fact that $R_n$ is an
equivalence relation on $U_n$, rather than on $X$, does not prevent us from
introducing the invariance notion expressed in \ref {PartialInvar}.  An
example of a function obeying this property is given by $\iota _n(f)$, where
$f$ is any function in $\Mp (U_n)$ which is constant on each
$R_n$-equivalence class.

Some elementary properties of these extended notions are in order.

\state Proposition \label UsefulPropTwo
  Given $n\in {\bf N}$, one has for all $f, g\in \MpX $, that
  \izitem
  \zitem $F_n\big (\iota _n(\rho _n)\big ) = \iota _n(\zeta _n)$,
  \zitemlbl QnOne $Q_n(1) = \charac {Y_n}$,
  \zitem $F_n(f) = F_n(\charac {U_n}f)$,
  \zitemlbl QnYn $Q_n(f) = Q_n(\charac {U_n}f) = Q_n(\charac {Y_n}f)$,
  \zitemlbl QnYnTwo $Q_n(f) = \charac {Y_n}Q_n(f)$,
  \zitemlbl FnQnNInvar $F_n(f)$ and $Q_n(f)$ are $R_n$-invariant and vanish
off $U_n$,
  \zitemlbl RkInvar if $f$ is $R_n$-invariant and vanishes off $U_n$, then
$f$ is $R_k$-invariant for every $k\leq n$,
  \zitemlbl CondexpForFQ if $g$ is $R_n$-invariant, then $F_n (gf)=gF_n
(f)$, and $Q_n (gf)=gQ_n (f)$,
  \zitem if $A$ is an $R_n$-invariant subset of $U_n$, then, as a subset of
$X$, $A$ is $R_k$-invariant for every $k\leq n$,
  \zitem $F_n$ and $Q_n$ are $\sigma $-additive.

\Proof Left for the reader.
  \endProof

\bigskip

Recalling that
  $$
  Z_n = Z_{\rho _n} = \{x\in U_n: \zeta _n(x)=\infty \},
  $$
  and that
  $$
  Y_n = U_n\sminus Z_{\rho _n} = \{x\in U_n: \zeta _n(x)<\infty \},
  $$
  we will now study certain relations among these sets, and we begin with
the following auxiliary result.

\state Lemma \label Decompos
  For every $0\leq n\leq m$, and for every $x$ in $U_m$, there exists a
subset $\Lambda \subseteq R_m(x)$, such that $x\in \Lambda $, and
  $$
  R_m(x) = \bigsqcup _{\lambda \in \Lambda } R_n(\lambda ),
  $$
  the square cup denoting disjoint union.

\Proof We first claim that if $C_n$ and $C_m$ are equivalence classes for
$R_n$ and $R_m$, respectively, then
  $$
  C_n\cap C_m\neq \emptyset \IMPLY C_n\subseteq C_m.
  $$ To see this, choose $z\in C_n\cap C_m$, and let $y\in C_n$.
  Then
  $$
  (y,z)\in R_n\cap (U_n\times U_m)\explain {Crazy}\subseteq R_m,
  $$
  so $y\in C_m$.
  This said, we see that the $R_m$-equivalence class of $x$ splits as the
union of $R_n$-equivalence classes, whence the conclusion.
  \endProof

The promissed relations among the $Z_n$ and the $Y_n$ are in order.

\state Proposition \label InclusionY
  \iaitem
  \aitem
  For every $m\geq 1$, and every $x$ in $U_m$, one has that \quad
  $
  e^{k_m(x)} \zeta _{m-1}(x)\leq \zeta _m(x).
  $
  \aitem If\/ $0\leq n\leq m$, then $Z_n\cap U_m\subseteq Z_m$,
  \aitem If\/ $0\leq n\leq m$, then $Y_m\subseteq Y_n$.

\Proof In order to prove (a) write
  $$
  R_m(x) = \bigsqcup _{\lambda \in \Lambda } R_{m-1}(\lambda ),
  $$
  where $x\in \Lambda \subseteq R_m(x)$, by \ref {Decompos}.  So
  $$
  \zeta _m(x) = \sum _{y\in R_m(x)} e^{h_m(y)} =
  \sum _{\lambda \in \Lambda }\ \sum _{y\in R_{m-1}(\lambda )} e^{h_m(y)} =
  \sum _{\lambda \in \Lambda }\ \sum _{y\in R_{m-1}(\lambda )}
e^{h_{m-1}(y)}e^{k_m(y)} = \cdots
  $$

For every $\lambda \in \Lambda $, and $y\in R_{m-1}(\lambda )$, notice that
  $$
  (y,\lambda )\in R_{m-1}\cap (U_m\times U_m)\explain {PotInv}\imply
k_m(y)=k_m(\lambda ),
  $$
  so the above equals
  $$
  \cdots =
  \sum _{\lambda \in \Lambda } e^{k_m(\lambda )} \kern -10pt \sum _{y\in
R_{m-1}(\lambda )} e^{h_{m-1}(y)} =
  \sum _{\lambda \in \Lambda } e^{k_m(\lambda )} \zeta _{m-1}(\lambda ) \geq
  e^{k_m(x)} \zeta _{m-1}(x).
  $$
  This proves (a).
  In order to prove (b), observe that under the hypothesis of (a) we have
that
  $$
  \zeta _{m-1}(x)=\infty \IMPLY \zeta _m(x)=\infty ,
  $$
  from where we trivially deduce that $Z_{m-1}\cap U_m\subseteq Z_m$.
Assuming now that $0\leq n\leq m$, we will prove (b) by induction on $m-n$.

In order to do this, notice that the case ``$m-n=0$'' is immediate, while
the case ``$m-n=1$'' has just been proved.  When $m-n>1$, we then have that
  $$
  Z_n\cap U_m\subseteq Z_n\cap U_{m-1}\cap U_m \subseteq Z_{m-1}\cap U_m
\subseteq Z_m,
  $$
  taking care of (b).  With respect to (c), we have
  $$
  Y_m = U_m\sminus Z_m = U_m \cap Z_m^c \explica {(a)} \subseteq
  U_m \cap (Z_n\cap U_m)^c \subseteq
  U_m \cap (Z_n^c\cup U_m^c) \quebra =
  (U_m \cap Z_n^c)\cup (U_m \cap U_m^c) =
  U_m \cap Z_n^c \subseteq
  U_n \cap Z_n^c =
  U_n\sminus Z_n= Y_n.
  \closeProof
  $$
  \endProof

There are many situations in the present context in which not necessarily
increasing sequences satisfy some \emph {increasing-like} property as we
look inside the appropriate $U_m$.  For example, when $n\leq m$, there is no
comparisson between $R_n$ and $R_m$, as sets, but when we restrict $R_n$ to
$U_m$, that is, when we consider $R_n\cap (U_m\times U_m)$, we get a subset
of $R_m$.  Similarly there is no comparission between $Z_n$ and $Z_m$, but
as seen above, $Z_n\cap U_m\subseteq Z_m$.

\medskip We next present some crucial properties of the $F_n$ and the $Q_n$.

\state Proposition \label LessElementary
  If\/ $0\leq n\leq m$, and if $f, g\in \MpX $, then
  \izitem
  \zitem $F_m\big (fF_n(g)\big )= F_m\big (F_n(f)g\big ), $
  \zitem $Q_m\big (Q_n(f)\big ) = Q_m(f) = Q_n\big (Q_m(f)\big ).$

\Proof Addressing (i), since both sides vanish on $X\sminus U_m$, by
definition, it is enough to prove that they agree on $U_m$.  Given $x$ in
$U_m$, we have
  $$
  F_m\big (fF_n(g)\big )\calcat x =
  \sum _{y\in R_m(x)}f(y)\sum _{z\in R_n(y)}g(z) = \cdots
  \equationmark FmFnfg
  $$
  Enploying \ref {Decompos} we write
  $$
  R_m(x) = \bigsqcup _{\lambda \in \Lambda } R_n(\lambda ),
  $$
  where $\Lambda \subseteq R_m(x)$, so \ref {FmFnfg} equals
  $$
  \cdots =
  \sum _{\lambda \in \Lambda }\ \sum _{y\in R_n(\lambda )}f(y)\sum _{z\in
R_n(y)}g(z) =
  \sum _{\lambda \in \Lambda }\ \sum _{y\in R_n(\lambda )}\ \sum _{z\in
R_n(\lambda )}f(y)g(z) \quebra =
  \sum _{\lambda \in \Lambda }\ \sum _{z\in R_n(\lambda )} g(z)\sum _{y\in
R_n(z)}f(y) =
  \sum _{z\in R_m(x)} g(z)\sum _{y\in R_n(z)}f(y) =
  F_m\big (gF_n(f)\big )\calcat x.
  $$
  This proves (i).  In order to prove (ii), recall that $h_m = \sum _{i=1}^m
k_i\restr {U_m}$.  So, defining
  $$
  \ell = \sum _{i=n+1}^m k_i\restr {U_m},
  $$
  we then have that $h_m = h_n + \ell $.  We next claim that
  $$
  (x,y)\in R_n\cap (U_m\times U_m) \IMPLY \ell (x)=\ell (y).
  $$
  This is because, for every $i=n+1,\ldots ,m$, we have that
  $$
  (x,y) \in
  R_n\cap (U_m\times U_m) \subseteq
  R_n\cap (U_n\times U_{i-1}) \subseteq R_{i-1},
  $$
  so $(x,y)$ also lies in $R_{i-1}\cap (U_i\times U_i)$, and this implies
that $k_i(x)=k_i(y)$, according to \ref {PotInv}.  In other words, $\ell $,
or rather $\iota _m(\ell )$, is $R_n$-invariant.

To prove (ii) we start with its left-hand-side, taking full advantage of the
abuse of language announced in \ref {ExtendedDefs.i}:
  $$
  Q_m\big (Q_n(f)\big ) \={CondexpForFQ}
  F_m\big (F_n(f\rho _n\zi _n)\rho _m\big )\zi _m \explica {(i)}=
  F_m\big (f\rho _n\zi _n F_n(e^\ell \rho _n)\big )\zi _m \={CondexpForFQ}
$$$$ =
  F_m\big (f\rho _n\zi _n e^\ell F_n(\rho _n)\big )\zi _m =
  F_m\big (f\rho _m\zi _n \zeta _n\big )\zi _m \={PropYr}
  F_m(f\charac {Y_n}\rho _m)\zi _m \quebra =
  Q_m(f\charac {Y_n}) \={QnYn}
  Q_m(f\charac {Y_m}\charac {Y_n}) \={InclusionY}
  Q_m(f\charac {Y_m}) = Q_m(f).
  $$
  With respect to the second equality in (ii), we have
  $$
  Q_n\big (Q_m(f)\big ) = Q_n\big (1 Q_m(f)\big )
\={UsefulPropTwo.\LocalFnQnNInvar ,\LocalRkInvar \&\LocalCondexpForFQ \stake
{5pt}}
  Q_n(1) Q_m(f) \={QnOne}
  \charac {Y_n} Q_m(f) \={QnYnTwo} $$$$=
  \charac {Y_n} \charac {Y_m} Q_m(f) \={InclusionY}
  \charac {Y_m} Q_m(f) = Q_m(f).
  \closeProof
  $$
  \endProof

We next present some useful properties of the dual operators $Q^*_n$.

\state Proposition \label LessElementaryStar
  Given $n$ in ${\bf N}$, and given any finite measure $\mu $ on $X$, one
has that
  \izitem
  \zitemlbl CondExpForStar $Q^*_n(g\mu ) = gQ^*_n(\mu )$, for every
$R_n$-invariant function $g$ in $\MpX $,
  \zitemlbl QnStarYn $Q_n^*(\mu ) = Q_n^*(\charac {Y_n}\mu )= \charac
{Y_n}Q_n^*(\mu )$,
  \zitemlbl QnStarAndPnStar $Q_n^*(\mu )\restr {U_n}=P_{\rho _n}^*(\mu
\restr {U_n})$,
  \zitemlbl QnOnInvarSet if $A$ is an $R_n$-invariant, Borel subset of $X$,
then $Q^*_n(\mu )(A)= \mu (A\cap Y_n)$,
  \zitemlbl Martingale if $m$ is another integer with $n\leq m$, then
$Q^*_m\big (Q^*_n(\mu )\big ) = Q^*_m(\mu ) = Q^*_n\big (Q^*_m(\mu )\big )$.

\Proof The first point follows easily from \ref {CondexpForFQ}.

\medskip \noindent \ref {LocalQnStarYn}: Pick any $f$ in $\MpX $.  Then
  $$
  \int _X f \d Q_n^*(\mu ) =
  \int _X Q_n(f) \d \mu \={QnYnTwo}
  \int _X Q_n(f)\charac {Y_n} \d \mu \quebra =
  \int _X Q_n(f)\d \charac {Y_n} \mu =
  \int _X f\d Q^*_n(\charac {Y_n} \mu ),
  $$
  proving the first identity in \ref {LocalQnStarYn}.  As for the second
one, we have
  $$
  \int _X f \d Q_n^*(\mu ) = \int _X Q_n(f) \d \mu \={QnYn}
  \int _X Q_n(\charac {Y_n}f) \d \mu =
  \int _X f\charac {Y_n} \d Q^*_n(\mu ),
  $$
  taking care of \ref {LocalQnStarYn}.

\medskip \noindent \ref {LocalQnStarAndPnStar}: Given $f\in \Mp (U_n)$, we
have
  $$
  \int _{U_n} f\d P_{\rho _n}^*(\mu \restr {U_n}) =
  \int _{U_n} P_{\rho _n}(f)\d \mu \restr {U_n} =
  \int _X \iota _n\big (P_{\rho _n}(f)\big )\d \mu \quebra =
  \int _X \iota _n\big (P_{\rho _n}(\iota _n(f)\restr {U_n})\big )\d \mu
\={RemarksOnExtendedDefs.b}
  \int _X Q_n\big (\iota _n(f)\big )\d \mu \quebra =
  \int _X \iota _n(f)\d Q_n^*(\mu ) =
  \int _{U_n} f\d Q_n^*(\mu )\restr {U_n},
  $$

\medskip \noindent \ref {LocalQnOnInvarSet}: We have
  $$
  Q^*_n(\mu )(A)=
  \int _X \charac A \d Q^*_n(\mu ) =
  \int _X Q_n(\charac A)\d \mu \Explain {CondexpForFQ}=
  \int _X \charac AQ_n(1)\d \mu \={QnOne}
  \int _X \charac A\charac {Y_n}\d \mu =
  \mu (A\cap Y_n).
  $$

\medskip \noindent \ref {LocalMartingale}: This is a direct consequence of
\ref {LessElementary.ii}.
  \endProof

We may now state an important quasi-invariance condition for measures on
$X$.

\state Theorem \label MainForQ
  Let $\mu $ be a finite measure on $X$.  Then $\mu $ is \dqi / if and only
if, for every $n\in {\bf N}$, one has that
  $$
  Q_n^*(\mu )=\charac {U_n}\mu .
  $$

\Proof We have already seen in \ref {ReduceToEtale} that $\mu $ is \dqi / if
and only each if $\mu \restr {U_n}$ is $D_n$-\qi /.  By \ref
{MainResultWithP} this is in turn equivalent to saying that
  $
  P_{\rho _n}^*(\mu \restr {U_n})=\mu \restr {U_n},
  $
  but in view of \ref {QnStarAndPnStar}, this is now the same as
  $$
  Q_n^*(\mu )\restr {U_n}=\mu \restr {U_n}.
  $$
  Since we know that $Q_n^*(\mu )$ vanishes on $X\sminus U_n$ by \ref
{QnStarYn}, the proof is concluded.
  \endProof

\newsection Existence of \qi / measures

\sectiontitle

\label ExistenceSection Having characterized \qi / measures in a concise way
in \ref {MainForQ}, we now discuss their existence.  This is a multi faceted
question manifesting itself in different ways on different parts of $X$.  It
is therefore convenient to break $X$ down into simpler pieces, so we shall
henceforth consider the following subsets
  $$
  V_n:=U_n\sminus U_{n+1}, \for n\in {\bf N},
  $$
  $$
  V_\infty :=\medcap _{n\in {\bf N}}U_n, \quad Z:=\medcup _{n\in {\bf
N}}Z_n,
  $$
  $$
  W_n:= V_n\sminus Z, \for n\in {\bf N}\cup \{\infty \},
  $$
  which we represent in the following diagram.  Please note that each $Z_n$
should be thought of as the largest shaded rectangle possessing the
indicated lower-left-hand corner.

\def \vline #1{\plot #1 0 #1 1000 /} \def \hline #1#2{\plot #1 #2 1500 #2 /}
\def \V #1#2#3{\put {$V_{#1}$} at #2 1050 \vline {#3}} \def \Z #1#2#3{\put
{$Z_{#1}\nearrow $} <-15pt, -7pt> at #2 #3} \def \W #1#2#3{\put {$W_{#1}$}
at #2 #3} \def \shade #1#2{\hshade #2 #1 1500 <0pt, 0pt, 0pt, 0pt> 0 #1 1500
/}

\beginmypix -100, 1600, -100, 1200 /

\setshadegrid span <1.2pt>

\def \xa {400} \def \xb {700} \def \xc {900} \def \xd {1000} \def \xe {1050}
\def \xf {1070} \def \xg {1080}

\def \ya {200} \def \yb {400} \def \yc {500} \def \yd {550} \def \ye {575}
\def \yf {587}

\plot 0 0 1500 0 1500 1000 0 1000 0 0 /

\V 0{200}{\xa } \V 1{550}{\xb } \V 2{800}{\xc } \V 3{950}{\xd }

\vline {\xe } \vline {\xf } \vline {\xg }

\put {$V_\infty $} at 1300 1050

\hline {\xa }{\ya } \hline {\xb }{\yb } \hline {\xc }{\yc } \hline {\xd
}{\yd } \hline {\xe }{\ye } \hline {\xf }{\yf }

\Z 1{\xa }{\ya } \Z 2{\xb }{\yb } \Z 3{\xc }{\yc }

\W 0{200}{500} \W 1{550}{600} \W 2{800}{700} \W 3{950}{800} \W \infty
{1300}{800}

\shade {\xa } {\ya } \shade {\xb } {\yb } \shade {\xc } {\yc } \shade {\xd }
{\yd } \shade {\xe } {\ye } \shade {\xf } {\yf }

\put {\eightrm Diagram \advseqnumbering (\current )} at 750 1130 \endmypix

The all important sets $U_n$'s may then be described as
  $$
  U_n = \bigsqcup _{n\leq m\leq \infty } V_m, \for n\in {\bf N},
  $$
  the square cup denoting disjoint union.
  In particular
  $$
  X = U_0 = \infunion V_n,
  $$
  whence also
  $$
  X = Z \ \sqcup \kern -3pt\infunion W_n.
  \equationmark DecompositionXWZ
  $$

\state Proposition \label InvariantPartition
  The following sets are $R$-invariant:
  \izitem
  \zitem $U_n$, for all $n\in {\bf N}$,
  \zitem $V_n$ and $W_n$, for all $n\in \Ninf $,
  \zitem $Z$.

\Proof Let us prove first that $U_m$ is $R$-invariant, for every $m$.  For
this suppose that $(x,y)\in R$, and $y\in U_m$, so we may pick some $n$ such
that $(x,y)\in R_n$.  Assuming initially that $n\leq m$, we have that
  $$
  (x,y)\in R_n\cap (U_n\times U_m)\explain {Crazy}\subseteq R_m\subseteq
U_m\times U_m,
  $$
  proving that $x\in U_m$, as desired.  Assuming now that $m\leq n$, the
conclusion comes even easier because
  $$
  (x,y)\in R_n\subseteq U_n\times U_n\subseteq U_m\times U_m,
  $$
  so again we have that $x\in U_m$.

Let us next prove that $Z$ is $R$-invariant.  So we pick $(x,y)\in R$, with
$y\in Z$, whence there are $n$ and $m$ such that $(x,y)\in R_n$, and $y\in
Z_m$.  Assuming initially that $n\leq m$, we have that
  $$
  (x,y)\in R_n\cap (U_n\times U_m)\explain {Crazy}\subseteq R_m.
  $$
  Since $\zeta _m$ is known to be $R_m$-invariant on $U_m$, it follows that
  $
  \zeta _m(x)=\zeta _m(y)=\infty ,
  $
  whence $x\in Z_m$, as required.

Assuming now that $m\leq n$, notice that
  $$
  y\in Z_m\cap U_n\explain {InclusionY.b}\subseteq Z_n,
  $$
  so the $R_n$-invariance of $\zeta _n$ implies that
  $$
  \infty =\zeta _n(y)=\zeta _n(x),
  $$
  and we conclude that
  $x\in Z_n\subseteq Z$.

Since the $R$-invariant sets clearly form a complete Boolean algebra, the
remaining statements follow.
  \endProof

The following result further streamlines the various quasi-invariance
conditions and it will be instrumental in the study of the existence
question.

\state Theorem \label CharacDLR
  Let $\mu $ be a finite measure on $X$, and for every $k\in \Ninf $, set
$\mu _k=\charac {W_k}\, \mu $.  Then $\mu $ is \dqi / if and only if all of
the following conditions hold:
  \izitem
  \zitem $\mu (Z)=0$,
  \zitem $Q^*_k(\mu _k) = \mu _k$, for every $k\geq 1$,
  \zitem $Q^*_i(\mu _\infty ) = \mu _\infty $, for every $i\geq 1$.

\Proof Assuming that $\mu $ is \dqi /, we have by \ref {ReduceToEtale} that
$\mu \restr {U_k}$ is $D_k$-\qi /, for every $k\in {\bf N}$, whence we
deduce from \ref {MainResult} that
  $$
  0 = \mu \restr {U_k}(Z_k) = \mu (Z_k),
  $$
  from where (i) follows.

Given $k\geq 1$, recall that $W_k$ is $R$-invariant, so $\mu _k$ is \dqi /
by \ref {InvarSetAndMeasure}.  In then follows from \ref {MainForQ}, that\fn
{Of course $Q^*_n(\mu _k)=\charac {U_n}\mu _k$, for all $n$, but so far we
only need the case $n=k$.}
  $$
  Q^*_k(\mu _k)=\charac {U_k}\mu _k = \charac {U_k}\charac {W_k}\mu =
\charac {W_k}\mu = \mu _k,
  $$
  proving (ii).

By the reasoning in the first sentence of the paragraph above, we also have
that $\mu _\infty $ is \dqi /.  So, again by \ref {MainForQ}, we have for
all $k\in {\bf N}$, that
  $$
  Q^*_k(\mu _\infty )=\charac {U_k}\mu _\infty = \charac {U_k}\charac
{W_\infty }\mu = \charac {W_\infty }\mu = \mu _\infty ,
  $$
  whence (iii).

Conversely, assuming that $\mu $ satisfies (i--iii), we will initially prove
that $\mu _k$ is \dqi / for all $k\in \Ninf $, via the characterization
provided by \ref {MainForQ}.  We must therefore prove that
  $$
  Q^*_n(\mu _k)=\charac {U_n}\mu _k, \for n\in {\bf N}.
  \equationmark CheckMainForQ
  $$
  When $k=\infty $, this is provided for by (iii), and the fact that
$W_\infty \subseteq U_n$, so it remains to prove \ref {CheckMainForQ} for
$k\in {\bf N}$.  Assuming first that $k<n$, we have that
  $$
  Q^*_n(\mu _k) \={LessElementaryStar.ii}
  Q^*_n(\charac {Y_n}\mu _k) =
  Q^*_n(\charac {Y_n}\charac {W_k}\mu _k) = 0,
  $$
  because $Y_n\cap W_k\subseteq U_n\cap W_k=\emptyset $, and likewise
$\charac {U_n}\mu _k=0$.  Assuming now that $n\leq k$, we have
  $$
  Q^*_n(\mu _k) \explica {(ii)}=
  Q^*_n\big (Q^*_k(\mu _k)\big ) \={LessElementaryStar.i}
  Q^*_k(\mu _k) \explica {(ii)}=
  \mu _k = \charac {U_k} \mu _k.
  $$

Note that the above use of (ii) is not quite correct because it has only
been assumed for $k\geq 1$.  However, when $k=0$, that condition holds
trivially because $Q^*_0$ is the identity transformation.
  This concludes the verification of \ref {CheckMainForQ}, hence showing
that $\mu _k$ is \dqi /.

Employing \ref {DecompositionXWZ} we have that
  $$
  \mu = \charac {Z} \mu + \medsum _{k\in \Ninf } \charac {W_k} \mu \
\explica {(i)}=\medsum _{k\in \Ninf } \mu _k,
  $$
  which is seen to be a \dqi / measure since each factor has this property,
which in turn is clearly preserved under sums.
  \endProof

Observe that any measure $\mu $ on $X$ which assigns zero mass to $Z$
satisfies
  $$
  \mu = \medsum _{k\in \Ninf } \mu _k,
  \equationmark SplitMu
  $$
  where each measure $\mu _k$ lives in $W_k$ (meaning that $\mu _k(X\sminus
W_k)=0$), namely $\mu _k=\charac {W_k}\mu $.  Conversely, if we are given a
collection $\{\mu _k\}_{k\in \Ninf }$ of measures on $X$, such that $\mu _k$
lives in $W_k$, then \ref {SplitMu} may be used to define a measure $\mu $
on $X$ which assigns measure zero to $Z$.  In other words there is a
one-to-one correspondence between the $\mu $'s and the collections $\{\mu
_k\}_{k\in \Ninf }$.  This said, observe that the conditions characterizing
a \dqi / measure in \ref {CharacDLR} consist of independent conditions on
each ``coordinate'' $\mu _k$.

In particular, if we fix any $k$ in $\Ninf $, and if we pick any measure
$\mu _k$ living on $W_k$, and satisfying the corresponding condition, namely
  \iaitem
  \aitem condition \ref {CharacDLR.ii}, in case $k\geq 1$,
  \aitem condition \ref {CharacDLR.iii}, in case $k=\infty $, or
  \aitem no condition at all, when $k=0$,
  \medskip \noindent then $\mu _k$, itself, is a \dqi / measure.  En
passant, we note that any measure living in $W_0$ is automatically \dqi /.
In any case, the existence question for \qi / measures should be split into
separate questions regarding the existence of \qi / measures living in each
$W_k$.  In case $k$ is finite, this question has a very simple answer:

\state Proposition
  Given an integer $n$, with $0\leq n<\infty $, there exists a \dqi /
probability measure living in $W_n$ if and only if\/ $W_n$ is nonempty.

\Proof
  Ignoring the blatantly obvious ``only if'' part, we deal only with the
``if'' part.  Assuming that $W_n$ is nonempty, choose any probability
measure $\nu $ living in $W_n$, e.g. a Dirac measure based on any point
chosen in $W_n$.  Setting $\mu =Q^*_n(\nu )$, notice that
  $$
  \mu (X) =
  \int _X 1 \d Q^*_n(\nu ) =
  \int _X Q_n(1) \d \nu \={QnOne}
  \int _X \charac {Y_n} \d \nu =
  \nu (Y_n) = 1,
  $$
  where the last equality is a consequence of the fact that
  $$
  W_n =
  V_n \sminus Z \subseteq
  U_n\sminus Z_n =
  Y_n.
  $$
  This shows that $\mu $ is a probability measure.

Since $W_n$ is $R_n$-invariant by \ref {InvariantPartition}, and since $\nu
$ lives in $W_n$, then $\mu $ also lives in $W_n$ by \ref {QnOnInvarSet}.
  For that reason we have that $\mu _k:= \charac {W_k}\mu = \delta _{nk}\mu
$, for every $k$ in ${\bf N}$.  So, in order to prove that $\mu $ is \dqi /,
we need only check \ref {CharacDLR.ii} for $k=n$, given that all of the
other conditions hold trivially.  To do this we compute
  $$
  Q^*_n(\mu ) =
  Q^*_n\big (Q^*_n(\nu )\big ) \= {LessElementaryStar.i}
  Q^*_n(\nu ) =
  \mu ,
  $$
  proving that $\mu $ is \dqi /.  \endProof

The existence question for \qi / measures living in $W_\infty $ is much more
subtle, not least because \ref {CharacDLR.iii} involves infinitely many
conditions.  The following result might have excessively rigid hypotheses,
but it is the best existence result we can offer in this generality:

\state Theorem
  Suppose that $W_\infty $ contains a compact, $R$-invariant, nonempty
subset $K$.  Suppose also that $\zeta _n\inv \restr K$ is
  continuous\fn {Please note that by saying that $\zeta _n\inv \restr K$ is
continuous we mean that it belongs to $C(K)$.  This should not be confused
with the much more stringent requirement that $\zeta _n\inv $ be continuous
at all points of $K$.}
  for every $n\in {\bf N}$.  Then there exists a \dqi / probability measure
living in $K$.

\Proof We begin by observing that $W_\infty \subseteq Y_n$, for every $n\in
{\bf N}$, because
  $$
  W_\infty =
  V_\infty \sminus Z \subseteq
  U_n\sminus Z_n =
  Y_n.
  $$

Denote by $P(X,K)$ the set of all probability measures on $X$ living in $K$.
This is clearly a nonempty set given that it contains any Dirac measure
$\delta _{x_0}$ with $x_0\in K$.  We claim that, for every $n$ in ${\bf N}$,
one has that
  $$
  Q_n^*\big (P(X,K)\big ) \subseteq P(X,K).
  $$
  In fact, if $\nu $ is in $P(X,K)$, then
  $$
  Q_n^*(\nu )(X) \={QnOnInvarSet} \nu (Y_n) = \nu (Y_n\cap K) = \nu (K) = 1,
  $$
  so we see that $Q_n^*(\nu )$ is a probability measure.  Moreover
  $$
  Q_n^*(\nu )(X\sminus K) \={QnOnInvarSet} \nu \big ((X\sminus K)\cap
Y_n\big ) \leq \nu (X\sminus K) = 0,
  $$
  so $Q_n^*(\nu )$ lives in $K$.
  Identifying $P(X,K)$ with the set $P(K)$ of all probability measures on
$K$ via the correspondence
  $$
  \mu \in P(X,K) \mapsto \mu \restr K\in P(K),
  $$
  we claim that, for every $g$ in $C(K)$, the function
  $$
  \mu \in P(X, K) \mapsto \int _K g \d Q^*_n(\mu ) \in {\bf C}
  $$
  is continuous on $P(K)$ relative to the topology induced by the weak*
topology of the dual of $C(K)$.

To prove it we first use Tietze's extension Theorem to produce a continuous
function $f$, defined on the whole of $X$, and whose restriction to $X$
coincides with $g$.  We further use Uhrysohn's Lemma to obtain a continuous
function $\varphi $ on $X$ such that $\varphi \restr K =1$, and whose
support is compact and contained in $U_n$.  Replacing $f$ by $\varphi f$ we
may then assume that the support of our originally chosen $f$ is compact and
contained in $U_n$.  We then have
  $$
  \int _K g \d Q^*_n(\mu ) =
  \int _X f \d Q^*_n(\mu ) =
  \int _X Q_n(f) \d \mu =
  \int _X i_n\big (P_{\rho _n}(f\restr {U_n})\big ) \d \mu \quebra =
  \int _{U_n} P_{\rho _n}(f\restr {U_n}) \d \mu \restr {U_n} =
  \int _{U_n} E_n(f\restr {U_n}\rho _n\zeta _n\inv ) \d \mu \restr {U_n}
\Explain {ERCondExp}=
  \int _{U_n} E_n(f\restr {U_n}\rho _n)\zeta _n\inv \d \mu \restr {U_n} =
  \int _K E_n(f\restr {U_n}\rho _n)\zeta _n\inv \d \mu .
  $$
  Recaling that $E_n(f\restr {U_n}\rho _n)$ is continuous by \ref
{IntroECont}, and that $\zeta _n\inv $ is continuous on $K$ by hypothesis,
the claim follows.

We will next prove that the fixed points of $Q^*_n$ in $P(X,K)$ form a
closed subset. To see this let $\mu $ be any measure in $P(X,K)$.  By the
uniqueness part of the Riesz-Markov, to say that $Q^*_n(\mu )=\mu $ is to
say that
  $$
  \int _K g \d \mu = \int _K g \d Q_n^*(\mu ), \for g \in C(K).
  \equationmark FirstFixed
  $$

Viewed as functions of the variable $\mu $, both sides of the above
expression are now known to be weak*-continuous, so the set of solutions to
this system on equations, as $g$ ranges in $C(K)$, form a closed set, hence
proving that the set of fixed points of $Q^*_n$ in $P(X,K)$ is indeed
weak*-closed.

Choosing any $\nu $ in $P(X,K)$, we define
  $$
  \mu _n=Q^*_n(\nu ),
  $$
  for every $n$ in ${\bf N}$.  Using Alaoglu's Theorem we may then find a
limit point, say
  $$
  \mu _\infty \in P(X, K),
  $$
  for the sequence $\{\mu _n\}_n$.
  Given two integers $n$ and $m$, with $0\leq n\leq m$, observe that
  $$
  Q^*_n(\mu _m) = Q^*_n\big (Q^*_m(\nu )\big ) \={Martingale} Q^*_m(\nu ) =
\mu _m,
  $$
  so we see that all but finitely many $\mu _m$'s are fixed points for
$Q^*_n$, hence the same holds for $\mu _\infty $, thanks to the closedness
of the set fixed points just proved.

Observing that $\mu _\infty $ lives in $K$, and that $K\subseteq U_n$, for
every $n$, we then have that
  $$
  Q^*_n(\mu ) = \mu = \charac {U_n}\mu ,
  $$
  so $\mu _\infty $ is \dqi / by \ref {MainForQ}.  \endProof

\newsection Renault-Deaconu groupoids

\sectiontitle

\label ExamplesSection In this section we will describe an example of \gap /
coming from generalized Renault-Deaconu groupoids \cite {cuntzlike}.  This
is in fact the main motivation for introducing and studying \gap /s.

We will henceforth fix a space $X$ satisfying \ref {Standing} and we will
suppose we are given an open subset $U\subseteq X$, and a map
  $$
  \sigma :U\to X,
  $$
  which we will assume to be a local homeomorphism.

Let $U_0=X$, and for each $n\geq 1$, put
  $$
  U_n = \big \{x\in U: \sigma (x)\in U_{n-1}\big \}.
  $$
  It is then easy to see that $U_n$ is effectively the domain of $\sigma
^n$, and that the $U_n$ form a collection of open subsets of $X$ satisfying
\ref {DefGap.i}.

The generalized Renault-Deaconu groupoid $\G _\sigma $ associated to $\sigma
$ was defined by Renault in \cite {cuntzlike}, and it consists of all
triples $(x,n,y)$ in $X\times {\bf Z}\times X$, such that there exist
$k,l\in {\bf N}$, satisfying
  $n=k-l$,
  $x\in U_k$, $y\in U_l$, and $\sigma ^k(x)=\sigma ^l(y)$.

The multiplication of two elements
  $(x,n,y)$ and $(z,m,w)$ in $\G _\sigma $ is defined only when $y=z$, in
which case the product is set to be $(x,n+m,w)$.
  The topology of $\G _\sigma $ is generated by the collection of subsets
  $$
  W_{k, l,A,B}=\big \{(x,k-l,y): x\in A,\ y\in B,\ \sigma ^k(x)=\sigma
^l(y)\big \},
  $$
  for all $k,l\in {\bf N}$, and all open subsets $A\subseteq U_k$, and
$B\subseteq U_l$.
  With this structure $\G _\sigma $ becomes an \'etale groupoid and we refer
the reader to \cite {cuntzlike} for more on $\G _\sigma $.

Let us next consider, for each $n\in {\bf N}$, the subset $R_n$ of
$U_n\times U_n$ defined by
  $$
  R_n=\big \{(x,y)\in U_n\times U_n: \sigma ^n(x)=\sigma ^n(y)\big \}.
  $$
  Recalling that $\sigma $ is a local homeomorphism, it is clear that
$\sigma ^n$ is a local homeomorphism from $U_n$ to $X$, so $R_n$ is easily
seen to be a \prp / equivalence relation on $U_n$.

\state Proposition One has that
  $$
  {\cal R} = \Big (\{U_n\}_{n\in {\bf N}}, \{R_n\}_{n\in {\bf N}}\Big ),
  $$
  is a \gap / on $X$.

\Proof
  All we need at this point is to verify \ref {DefGap.iii}.  So, suppose
that $n\leq m$, and let $(x, y)\in R_n\cap (U_n\times U_m)$.  The first
conclusion to be drawn is that both $x$ and $y$ lie in $U_n$, and that
$\sigma ^n(x)=\sigma ^n(y)$.  Besides, $y$ lies in $U_m$, so $y$ is also in
the domain of $\sigma ^m$.  This implies in particular that $\sigma ^n(y)$
lies in the domain of $\sigma ^{m-n}$, so the same holds for $\sigma ^n(x)$.
Consequently
  $$
  \sigma ^m(y) = \sigma ^{m-n}\big (\sigma ^n(y)\big ) = \sigma ^{m-n}\big
(\sigma ^n(x)\big ) = \sigma ^m(x),
  $$
  thus showing that $(x,y)\in R_m$.  \endProof

The Renault-Deaconu groupoid admits a continuous cocycle
  $$
  \delta :(x,n,y)\in \G _\sigma \mapsto n\in {\bf Z},
  $$
  whose kernel is therefore an open subgroupoid, which is clearly isomorphic
and homeomorphic to the \gap / groupoid $R=\bigcup _n R_n$, via the
homeomorphism sending each $(x,y)$ in $R$ to $(x,0,y)$.

In order to speak of \qi / measures on $\G _\sigma $, we need to introduce a
real-valued cocycle.  Recall from \cite {PaperOne} that if $h:U\to {\bf R}$
is a continuous function, often thought of as a \emph {potential} function,
one may define
  $$
  b(x,n-m, y) =
  \medsum _{i=0}^{n-1} h\big (\sigma ^i(x)\big ) - \medsum _{i=0}^{m-1}
h\big (\sigma ^i(y)\big ),
  \equationmark RDCocycle
  $$
  whenever $\sigma ^n(x)=\sigma ^m(y)$, obtaining in this way a well defined
continuous cocycle
  on $\G _\sigma $,
  taking values in the additive group of real numbers.

The restriction of $c$ to the subgroupoid $R$ is then evidently a continuous
cocycle on $R$, but we would instead like to introduce it from a different
point of view, in line with \ref {GapCocycle}.
  For each $n\geq 1$, let us define
  $$
  k_n(x)=h\big (\sigma ^{n-1}(x)\big ), \for x\in U_n.
  $$
  Notice that, for $x$ in $U_n$, the largest integer $i$ for which the
expression $h\big (\sigma ^i(x)\big )$ is guaranteed to be well defined is
$i=n-1$.  This is because $U_n=\hbox {dom}(\sigma ^n)$, so $\sigma ^n(x)$ is
well defined, which in turn implies that $\sigma ^{n-1}(x)$ is in $U$, and
hence it does makes sense to apply $h$ to $\sigma ^{n-1}(x)$.  However there
is no reason for $\sigma ^n(x)$ to lie in $U$, so $h(\sigma ^n(x))$ may not
be well defined.

\state Proposition The collection $\{k_n\}_{n\geq 1}$ defined above is a
potential for $\cal R$, in the sense of \ref {DefinePotential}.

\Proof To verify \ref {PotInv}, let $n\geq 1$ and choose
  $(x,y)\in R_{n-1}\cap (U_n\times U_n)$.  Then $\sigma ^{n-1}(x)=\sigma
^{n-1}(y)$, whence
  $
  k_n(x)=k_n(y),
  $
  as needed.
  \endProof

The associated cocycle is then given on any $(x,y)\in R_n$, by
  $$
  c(x,y) =
  c_n(x,y) =
  h_n(x) - h_n(y) =
  \medsum _{i=1}^n k_n(x) - k_n(y) \quebra =
  \medsum _{i=1}^n h\big (\sigma ^{i-1}(x)\big ) - h\big (\sigma
^{i-1}(y)\big ) =
  \medsum _{i=0}^{n-1} h\big (\sigma ^i(x)\big ) - h\big (\sigma ^i(y)\big
),
  \equationmark DefineCFromh
  $$
  which happens to be the restriction to $R$ of the cocycle $b$ defined in
\ref {RDCocycle}.

Observe that both the unit space of $\G _\sigma $ and that of $R$ may be
naturally identifyied with $X$.  So a given finite measure $\mu $ on $X$ may
be tested for quasi-invariance either relatively to the cocycle $e^b$ on $\G
_\sigma $, or to the cocycle $e^c$ on $R$.

\definition \label DefineConformalAndDLR
  Let $\mu $ be a finite measure on $X$.  We will say that $\mu $ is
  \izitem
  \zitem a \emph {conformal} measure when it is \qi / relatively to the
cocycle $e^b$ on $\G _\sigma $,
  \zitem a \emph {DLR} measure when it is \qi / relatively to the cocycle
$e^c$ on $R$.

Since $R$ is a subgroupoid of $\G _\sigma $, and since $e^c$ is the
restriction of $e^b$ to $R$, it is immediate that:

\state Proposition
  Every conformal measure on $X$ is a DLR measure.

\newsection Eigenmeasures

\def \X {X} \def \Xs #1{\X }

\sectiontitle

As in the previous section we let $X$ be a space satisfying \ref {Standing},
$U\subseteq X$ be an open set, and $\sigma :U\to X$ be a local
homeomorphism.
  Our goal here is to show that every eigenmeasure for Ruelle's operator is
a DLR measure.

For each $f$ in $\Mp (U)$, and for every $x$ in $\X $, define
  $$
  L(f)\calcat x = \sum _{t\in \sigma \inv (x)} f(t),
  $$
  so that
  $L$ becomes a map
  $$
  L:\Mp (U)\to \Mp (\X ).
  $$

Regarding the expression defining $L(f)$ above, observe that if $x$ is not
in $\sigma (U) $, then $\sigma \inv (x)$ is the empty set, whence there are
no summands at all, hence the sum turns out to be zero.  In other words,
$L(f)$ vanishes outside $\sigma (U)$.

  We also consider the operator
  $$
  \alpha :\Mp (\X )\to \Mp (U),
  $$
  given by $\alpha (f) = f\circ \sigma $.
  We finally define
  $$
  E:\Mp (U)\to \Mp (U),
  $$
  by
  $$
  E(f)\calcat x = \sum _{\sigma (t)=\sigma (x)} f(t), \for x\in U.
  $$

For $f$ in $\Mp (U)$ and any $x$ in $U$, observe that,
  $$
  \alpha \big (L(f)\big )\calcat x =
  L(f)\calcat {\sigma (x)} =
  \sum _{t\in \sigma \inv (\sigma (x))} f(t) =
  \sum _{\sigma (t)=\sigma (x)} f(t) =
  E(f)\calcat x,
  $$
  so we see that
  $$
  \alpha \circ L=E.
  \equationmark ALE
  $$
  Anoter useful property is
  $$
  L(\alpha (g)f)=gL(f), \for g\in \Mp (\X ), \for f\in \Mp (U),
  \equationmark CondexpL
  $$
  which the reader will have no difficulty in checking.

We shall also fix a continuous function
  $$
  \rho :U\to {\bf R},
  $$ satisfying $\rho (x)>0$, for all $x$ in $U$.  Here $\rho $ is supposed
to play the role of $e^h$, where $h$ is the function used for creating the
cocycle $b$ in \ref {RDCocycle}.
  The operator
  $$
  L_\rho :\Mp (U)\to \Mp (\X ),
  $$
  defined by the formula
 $L_\rho (f)=L(\rho f)$, is then the analogue of Ruelle's operator in the
present situation.

\state Lemma \label OneStep
  Suppose that
  \izitem
  \zitem $\mu $ is a measure on $X$,
  \zitem $\lambda $ is a nonzero scalar,
  \medskip \noindent such that
  $$
  \int _\X L(\rho f)\d \mu = \lambda \int _U f\d \mu , \for f\in \Mp (U).
  \equationmark EigenMeasure
  $$
  Then
  $$
  \int _U E(\rho )f \d \mu = \int _U E (\rho f) \d \mu , \for f\in \Mp (U).
  \equationmark DLR
  $$

\Proof Given $g$ in $\Mp (\X )$, plug $f=\alpha (g)$ in \ref {EigenMeasure},
to get
  $$
  \int _U \alpha (g)\d \mu =
  \lambda \inv \int _\X L\big (\alpha (g)\rho \big )\d \mu \={CondexpL}
  \lambda \inv \int _\X gL(\rho )\d \mu .
  \equationmark WhatAlphaDoes
  $$ Working from the right-hand-side of \ref {DLR}, we have
  $$
  \int _U E (\rho f) \d \mu \={ALE}
  \int _U \alpha \big (L(\rho f)\big ) \d \mu \={WhatAlphaDoes}
  \lambda \inv \int _\X L(\rho f) L(\rho )\d \mu \Explain {CondexpL}=
  \lambda \inv \int _\X L\big (\rho f\alpha \big (L(\rho )\big )\big ) \d
\mu \={ALE}
  \lambda \inv \int _\X L\big (\rho fE(\rho )\big ) \d \mu \={EigenMeasure}
  \int _U fE(\rho )\d \mu .
  \closeProof
  $$
  \endProof

For each $n\geq 1$, let $U_n$ be the domain of $\sigma ^n$, so that the map
  $$
  \sigma ^n:U_n\to X,
  $$
  is a local homeomorphism so it could be used in place of $\sigma $ in
order to define all of the above ingredients.  To be precise these are:

\bigskip \Bitem $\ds
  L_n:\Mp (U_n)\to \Mp (\Xs n )$,\quad given by \quad $L_n(f)\calcat x = \ds
\sum _{t\in \sigma ^{-n}(x)} f(t),
  $ \Bitem
  $
  \stake {15pt}\pilar {10pt}
  \alpha _n:\Mp (\Xs n )\to \Mp (U_n)$,\quad given by \quad $\alpha _n(f) =
f\circ \sigma ^n,
  $ \Bitem
  $\ds
  E_n:\Mp (U_n)\to \Mp (U_n)$,\quad given by \quad $E_n(f)\calcat x = \ds
\sum _{\sigma ^n(t)=\sigma ^n(x)} f(t).
  $

\bigskip \noindent We shall also let
  $$
  \rho _n = \rho \alpha \big (\rho \big )\alpha _2(\rho )\cdots \alpha
_{n-1}(\rho ).
  $$

  Regarding the product defining $\rho _n$ above,
  observe that each $\alpha _i(\rho )$ is a member of $\Mp (U_{i+1})$, so
they all may be restricted to $U_n$ before the multiplication is performed,
resulting of course in a member of $\Mp (U_n)$.

\state Proposition \label AnulaFora
  For every $n, m\geq 0$, one has that
  \Zitem $L_n\big (\Mp (U_{n+m})\big )\subseteq \Mp (U_m)$, and
  \zitem $L_m\circ L_n = L_{n+m}$.

\Proof
  Given $f$ in $\Mp (U_{n+m})$, and $x$ in $X\sminus U_m$, we must prove
that
  $L_n(f)\calcat x = 0$.
  Arguing by contradiction, suppose this is not so.  Therefore there exists
at least one $t\in \sigma ^{-n}(x)$, such that $f(t)\neq 0$, so
  it follows that
  $$
  x = \sigma ^n(t) \in \sigma ^n(U_{n+m})\subseteq U_m,
  $$
  a contradiction, proving (i). In order to prove (ii), pick
  $f$ in $\Mp (U_{n+m})$, and $x$ in $X$.
  We then have
  $$
  L_m\big (L_n(f)\big )\calcat x =
  \sum _{t\in \sigma ^{-m}(x)} L_n(f)\calcat t =
  \sum _{t\in \sigma ^{-m}(x)} \sum _{s\in \sigma ^{-n}(t)} f(s) =
  \sum _{t\in \sigma ^{-n-m}(x)} f(s) = L_{n+m}(f)\calcat x.
  \closeProof
  $$
  \endProof

\state Lemma \label HalfStep
  Let $\mu $ be a measure on $X$ satisfying \ref {EigenMeasure}.  Then
  $$
  \int _{\Xs n }L_n(\rho _nf)\d \mu = \lambda ^n\int _{U_n} f\d \mu , \for
f\in \Mp (U_n).
  $$

\Proof For $f$ in $\Mp (U_n)$, we have by induction that
  $$
  \int _{\Xs n }L_n(\rho _nf)\d \mu =
  \int _{\Xs n }L_{n-1}\big (L(\rho \alpha (\rho _{n-1})f)\big )\d \mu =
  \int _{\Xs n }L_{n-1}\big (\rho _{n-1}L(\rho f)\big )\d \mu \quebra =
  \lambda ^{n-1}\int _{U_{n-1}} L(\rho f)\d \mu =
  \lambda ^{n-1}\int _{U} L(\rho f)\d \mu \={EigenMeasure}
  \lambda ^n\int _U f \d \mu =
  \lambda ^n\int _{U_n} f \d \mu .
  \closeProof
  $$ \endProof

\state Corollary \label CoroGeral
  Let $\mu $ be a measure on $X$ satisfying \ref {EigenMeasure}.  Then for
any $n$ in ${\bf N}$, one has that
  $$
  \int _{U_n} E_n(\rho _n)f \d \mu = \int _{U_n} E_n (\rho _nf) \d \mu ,
\for f\in \Mp (U_n).
  $$

\Proof
  Follows immediately by applying \ref {OneStep} to $\sigma ^n$ and $\rho
_n$.  \endProof

Given a continuous function $h:U\to {\bf R}$, let $\rho =e^h$, and recall
from \ref {DefineConformalAndDLR} that a finite measure $\mu $ on $X$ is
said to be a DLR measure for $h$ if it is \qi / relative to the cocycle
$e^c$ on $R$, where $c$ is given in terms of $h$ by \ref {DefineCFromh}.
The cocycle $e^c$ is in fact a common extension of the cocycles $e^{c_n}$
defined on each $R_n$ by
  $$
  e^{c_n(x,y)} =
  \exp \Big (\medsum _{i=0}^{n-1} h\big (\sigma ^i(x)\big ) - h\big (\sigma
^i(y)\big )\Big ) =
  {\rho (x) \rho (\sigma (x))\cdots \rho (\sigma ^{n-1}(x)) \over \rho (y)
\rho (\sigma (y))\cdots \rho (\sigma ^{n-1}(y))} = {\rho _n(x) \over \rho
_n(y)}.
  $$

We then see that, if $\mu $ is a finite measure satisfying the conclusions
of \ref {CoroGeral}, then $\mu \restr {U_n}$ satisfies \ref {Harmonic}
relative to $\rho _n$, so it is \qi / for $e^{c_n}$ by \ref {MainResult}.
Employing \ref {ReduceToEtale} it then follows that $\mu $ is $e^c$-\qi /,
hence a DLR measure.  Summarizing we obtain the following:

\state Theorem \label EigenDLR
  Let $X$ be a locally compact metric space, $U$ be an open subset of $X$,
and $\sigma :U\to X$ be a local homeomorphism.  Choosing any continuous
potential $h:U\to {\bf R}$, let $\rho =e^h$.  Then any finite measure on $X$
which is an eigenvalue for the corresponding Ruelle operator, meaning that
it satisfies \ref {EigenMeasure} with a nonzero eigenvalue $\lambda $, is
necessarily a DLR measure.

\Proof
  Follows immediately by applying \ref {CoroGeral} to $\sigma ^n$ and $\rho
_n$.  \endProof

It should be noted that Corollary \ref {CoroGeral} may be seen as a
generalization of Theorem \ref {EigenDLR} to measures which are not
necessarily finite, as long as we redefine the notion of DLR measures as
those which satisfy the conclusions of \ref {CoroGeral}.  However, since our
theory of DLR measures was developed only for finite measures, the various
equivalent conditions for a measure to be DLR have not been proved here for
infinite measures.

\def \sectiontitle \par {\stno = 0
    \goodbreak \bigbreak
    \noindent {\bf \AppendixId .\enspace \secname .}
    \nobreak \medskip \noindent }

\def \current {\AppendixId \ifnum \stno = 0 \else .\number \stno \fi }

\newsection Appendix: Elementary remarks about Measure Theory

\def \AppendixId {A}

\sectiontitle

In the final two sections of this work we make some elementary remarks,
mainly to fix our notation.
  By a measurable space we shall mean a pair $ (X, {\cal B}) $, consisting
of a nonempty set $X$, and a $\sigma $-algebra $\cal B$ of subsets of $X$.

\definition \label IntroMp
  Given a measurable space $ (X, {\cal B}) $, we shall denote by $\Mp (X,
{\cal B}) $ the set of all $\cal B$-measurable functions
  $$
  f:X\to [0,+\infty ].
  $$

\definition \label IntroSigmaAdd
  Given measurable spaces $ (X, {\cal B}) $ and $(Y, {\cal C}) $, a
positively homogeneous map
  $$
  T:\Mp (X, {\cal B}) \to \Mp (Y, {\cal C})
  $$
  is said to be \emph {$\sigma $-additive} if for any sequence $\{f_n\}_n$
in $\Mp (X, {\cal B}) $, we have that
  $$
  T \Big (\sum _{n=1}^\infty f_n \Big ) = \sum _{n=1}^\infty T (f_n),
  $$
  all sums being interpreted as pointwise sums.

\state Proposition \label Daniell
  Let $ (X, {\cal B}) $ and $(Y, {\cal C}) $ be measurable spaces, and let
  $$
  T:\Mp (X, {\cal B}) \to \Mp (Y, {\cal C})
  $$
  be a $\sigma $-additive map.  Then:
  \izitem
  \zitem
  If\/ $\{f_n\}_n$ is a non-decreasing sequence of functions in $\Mp (X,
{\cal B}) $, then
  $$
  T \big (\lim _{n\to \infty } f_n\big ) = \lim _{n\to \infty }T (f_n),
  $$
  all limits being interpreted as pointwise limits.  \zitemmark MonotoneForE
  \zitem For every measure\fn {All measures in this work are assumed to be
$\sigma $-additive and positive.} $\nu $ on $(Y, {\cal C}) $, there exists a
measure $T^*(\nu )$ on $ (X, {\cal B}) $, such that
  $$
  \int _Xf\d T^*(\nu ) = \int _YT(f)\d \nu , \for f\in \Mp (X, {\cal B}) .
  $$

\Proof Regarding (i), define $g_1=f_1$, and for each $n\geq 2$, define
  $$
  g_n=f_n(x)-f_{n-1}(x), \for x\in X.
  $$
  Observe that, in case $f_{n-1}(x)=\infty $, then necessarily also
$f_n(x)=\infty $, in which case we adopt the convention according to which
$g_n(x)=0$.  The functions $g_n$ so defined are then $\cal B$-measurable and
non-negative, and we have that $f_{n-1}+g_n=f_n$.  Consequently $f_n=\sum
_{i=1}^ng_i$, so
  $$
  T \big (\lim _{n\to \infty } f_n\big ) =
  T \Big (\medsum _{i=1}^\infty g_i\Big ) =
  \medsum _{i=1}^\infty T (g_i) =
  \lim _{n\to \infty }\medsum _{i=1}^nT (g_i) =
  \lim _{n\to \infty }T (f_n),
  $$
  proving (i).
  In order to prove (ii), for every $f$ in $\Mp (X, {\cal B}) $, define
  $$
  I(f) = \int _XT(f)\d \nu .
  $$

  We then claim that, given any sequence $\{f_n\}_n$ in $\Mp (X, {\cal B})
$, one has that
  $$
  I\big (\medsum _{n=1}^\infty f_n\big ) =
  \medsum _{n=1}^\infty I(f_n).
  $$
  This is a consequence of the $\sigma $-additivity of $T$ and the monotone
convergence Theorem, as follows:
  $$
  I\big (\medsum _{n=1}^\infty f_n\big ) =
  \int _XT \big (\medsum _{n=1}^\infty f_n\big ) \d \nu =
  \int _X \medsum _{n=1}^\infty T(f_n)\d \nu =
  \medsum _{n=1}^\infty \int _X T(f_n)\d \nu =
  \medsum _{n=1}^\infty I(f_n),
  $$
  thus proving the claim.
  Defining $\mu (E)=I(\charac E)$, for every $E$ in $\cal B$, it then
follows that $\mu $ is a $\sigma $-additive measure on $X$, and clearly
  $$
  I(f) = \int _Xf\d \mu ,
  \equationmark IntTau
  $$
  for every simple function $f$ in $\Mp (X, {\cal B}) $.

We then claim that \ref {IntTau} holds for every $f$ in $\Mp (X, {\cal B})
$.  To see this, recall that every such $f$ may be written as the pointwise
limit of a non-decreasing sequence $\{s_n\}_n$ of simple functions in $\Mp
(X, {\cal B}) $ \cite [Section 18.1]{RK}.  So
  $$
  I(f) =
  \int _XT\big (\lim _{n\to \infty } s_n\big )\d \nu \= {LocalMonotoneForE}
  \int _X\lim _{n\to \infty } T(s_n)\d \nu \explica {(*)}=
  \lim _{n\to \infty } \int _XT(s_n)\d \nu \quebra =
  \lim _{n\to \infty } I(s_n)=
  \lim _{n\to \infty } \int _Xs_n\d \mu =
  \int _Xf\d \mu .
  $$
  In time, we observe that $(^*)$, above, is justified by the monotone
convergence Theorem and the easily proven fact that $\{T(s_n)\}_n$ is a
non-decreasing sequence.

  This proves our claim and we then have for every $f$ in $\Mp (X, {\cal B})
$ that
  $$
  \int _Xf\d \mu = I(f) = \int _XT(f)\d \nu ,
  $$
  so it is enough to put $T^*(\nu ):=\mu $.
  \endProof

Before closing this section let us comment on two related notions which will
be used often.

\state Remark \label RestrictingMeasures
  \rm In this work we shall consider two similar, but inequivalent, ways of
restricting a measure on a space $X$ to a Borel subset $Y\subseteq X$.  The
first one, officially called the \emph {restriction of $\mu $ to $Y$},
consists of the measure denoted $\mu \restr Y$, defined on the $\sigma
$-algebra ${\cal B}(Y)$ of all Borel subsets of $Y$, by
  $$
  \mu \restr Y(A) = \mu (A), \for A\in {\cal B} (Y).
  $$ The second one, which we will denote by $\charac Y\mu $, is nothing but
the well known measure obtained by multiplying the measurable function
$\charac Y$ by the measure $\mu $. Recall that the domain of $\charac Y\mu $
is still ${\cal B}(X)$, as opposed to ${\cal B}(Y)$, and
  $$
  \charac Y\mu (A)=\int _A1_Y\d \mu = \mu (A\cap Y),
  $$
  for all $A\in {\cal B}(X)$.

\newsection Appendix: Elementary remarks about Measure Theory and Topology

\def \AppendixId {B}

\sectiontitle

\state Proposition
   Every space $X$ satisfying \ref {Standing} is $\sigma $-compact.

\Proof For every $x$ in $X$, let $U_x$ be a relatively compact, open
neighborhood of $x$.  By reducing $U_x$ a bit we may assume that it belongs
to some previously chosen countable base $B$ of open sets for $X$.  It then
follows that
  $$
  \{\overline U_x:x\in X\}
  $$
  is a family of compact sets covering $X$.  This family is countable (even
though it might be indexed on an uncountable set) because the $U_x$'s belong
to the countable base $B$.
  \endProof

The reason for restricting ourselves to \ref {Standing} is to simplify some
aspects of measure theory.  In this short section we will explain exactly
what we mean by this.

Recall that the \emph {Borel $\sigma $-algebra}, denoted ${\cal B}(X)$, is
the $\sigma $-algebra of subsets of $X$ generated by the closed subsets.  On
the other hand, the \emph {Baire $\sigma $-algebra} \cite [21.6]{RK},
denoted ${\cal B}a(X)$, is the smallest $\sigma $-algebra of subsets of $X$
for which the functions in $C_c (X)$ are measurable.

If one is interested in the measurability properties of none other that
continuous, compactly supported functions, the Baire $\sigma $-algebra is
the most appropriate one to be considered.
  The Baire $\sigma $-algebra is known to be generated by the compact
$G_\delta $ subsets of $X$ \cite [Theorem 21.21]{RK}.

The first advantage of working with \ref {Standing} is as follows:

\state Lemma {\rm (See also \cite [Theorem 21.20]{RK})} \
  Suppose that $X$ is as in \ref {Standing}.  Then the Baire and Borel
$\sigma $-algebras on $X$ coincide.

\Proof It is clear that ${\cal B}a(X)\subseteq {\cal B}(X)$, so we need only
wory about the reverse inclusion.  In order to do so it is enough to prove
that every closed set $F\subseteq X$ is Baire-measurable.  Temporarily
assuming that $F$ is moreover compact, pick any compatible metric $d$ on $X$
and define
  $$
  U_n=\{x\in X:d(x,F)<1/n\}.
  $$
  Each $U_n$ is then open, and $F=\bigcap _nU_n$, so we see that $F$ is a
compact $G_\delta $, hence Baire-measurable.

If $F$ is now any closed set, use the fact that $X$ is $\sigma $-compact to
choose a countable family $\{K_n\}_n$ of compact subsets of $X$ such that
$X=\bigcup _n K_n$.
  Then $F\cap K_n$ is a compact set, and hence Baire-measurable by the first
part of this proof.  Since $F = \bigcup _n F\cap K_n$, it follows that $F$
is Baire-measurable.
  \endProof

A measure $\mu $ defined on ${\cal B}(X)$ is caled a \emph {Borel measure}
when it assigns finite measure to every compact set \cite [Section
21.3]{RK}.  If $\mu $ is instead defined on ${\cal B}a(X)$, it is called a
\emph {Baire measure} provided it is finite on compact (Baire-measurable)
sets \cite [Section 21.6]{RK}.

One of the main applications of Borel measures in Analysis is the
Riesz-Markov Theorem \cite [Section 21.6]{RK} which states that each
positive linear functional on $C_c(X)$ is given by the integration against a
unique \emph {regular} Borel measure (a regular Borel measure is also called
a \emph {Radon} measure).

Another reason to work under \ref {Standing} is that, in this case, every
Baire measure is regular \cite [Theorem 21.27]{RK}.  Since we now know that
the Baire and Borel $\sigma $-algebras coincide, we deduce that:

\state Lemma \label AllRegular
  {\rm (See also \cite [Theorem 21.20]{RK})} \
  Under \ref {Standing}, every Borel measure on ${\cal B}(X)$ is regular.

\state Proposition \label ConditionBorel
  Let $\mu $ be any measure on ${\cal B}(X)$.  Then $\mu $ is a Borel
measure if and only if
  $$
  \int _Xf\d \mu <\infty ,
  $$
  for all $f$ in $C_c^+(X)$.

\Proof We verify only the ``if'' part.  For this, let $K$ be any compact
subset of $X$.  Using local compactness one may find a relatively compact,
open set $U$ such that $K\subseteq U$.  By Uryshon's Theorem let $f:X\to
[0,1]$ be a continuous function vanishing off $U$, and such that $f=1$ on
$K$.  If follows that $f\in C_c^+(X)$, whence
  $$
  \mu (K) = \int _X \charac K\d \mu \leq \int _X f\d \mu <\infty .
  $$
  Being finite on compact sets, $\mu $ is a Borel measure.
  \endProof

\state Proposition \label ApproxUnit
  There exists a sequence $\{\varphi _n\}_n$ of continuous, compactly
supported functions
  $$
  \varphi _n:X\to [0,1],
  $$
  such that $\varphi _n\leq \varphi _{n+1}$, for every $n$, and $\ds \lim
_{n\to \infty } \varphi _n=1$, pointwise.

\Proof Let $\{K_n\}_{n\in {\bf N}}$ be a sequence of compact subsets of $X$
such that
  $$
  K_n\subseteq \hbox {int}(K_{n+1}), \and X = \bigcup _nK_n.
  $$
  Using Uryhson, for each $n\in {\bf N}$, let
  $$
  \varphi _n:X\to [0,1]
  $$
  be a continuous function with $\varphi _n=1$ on $K_n$ and $\varphi _n=0$
on $X\setminus \hbox {int}(K_{n+1})$.  It is then easy to see that $\varphi
_n\leq \varphi _{n+1}$, and that $\{\varphi _n\}_n$ converges pointwise to 1
on $X$.
  \endProof

\references

\Bibitem PaperOne
  R. Bissacot, R. Exel, R. Frausino and T. Raszeja;
  Conformal measures on generalized Renault-Deaconu groupoids;
  preprint, 2018

\Article actions
  R. Exel;
  Inverse semigroups and combinatorial C*-algebras;
  Bull. Braz. Math. Soc. (N.S.), 39 (2008), 191-313

\Article infinoa
  R. Exel and M. Laca;
  Cuntz-Krieger algebras for infinite matrices;
  J. reine angew. Math., 512 (1999), 119-172

\Article eq
  R. Exel and A. Lopes;
  C*-algebras, approximately proper equivalence relations, and thermodynamic
formalism;
  Ergodic Theory Dynam. Systems, 24 (2004), 1051-1082

\Bibitem Renault
  J. Renault;
  A groupoid approach to C*-algebras;
  Lecture Notes in Mathematics vol.~793, Springer, 1980

\Bibitem cuntzlike
  J. Renault;
  Cuntz-like algebras;
  Proceedings of the 17th International Conference on Operator Theory
(Timisoara 98), The Theta Fondation, 2000

\Article ApRenault
  J. Renault;
  The Radon-Nikod\'ym problem for appoximately proper equivalence relations;
  Ergodic Theory Dynam. Systems, 25 (2005), no. 5, 1643-1672

\Bibitem RK
  H. L. Royden and P. M. Fitzpatrick;
  Real Analysis;
  fourth edition,
  Pearson Education Asia Ltd., 2010

\endgroup

\bye